\documentclass{gen-j-l}%
\usepackage{amssymb}
\usepackage{amsmath}
\usepackage{amsfonts}%
\usepackage{graphicx}
\newtheorem{theorem}{Theorem}[section]
\newtheorem{lemma}[theorem]{Lemma}
\theoremstyle{definition}
\newtheorem{definition}[theorem]{Definition}

\theoremstyle{remark}

\numberwithin{equation}{section}
\theoremstyle{plain}

\copyrightinfo{2001}{enter name of copyright holder}
\begin{document}
\title[Mixed Dimensional Compactness]{Mixed Dimensional Compactness with Dimension Collapsing from $\mathbb{S}%
^{n-1}$ Bundle Measures.}
\author{Simon P Morgan}
\address{University of Minnesota}
\email{morgan@math.umn.edu}
\urladdr{http://www.math.umn.edu/\symbol{126}morgan}

\begin{abstract}
We provide a measure based topology for certain unions of $C^{2}$ rectifiable
submanifolds of mixed dimensions in $\mathbb{R}^{n}$. In this topology lower
dimensional sets remain in the limit as measures when higher dimensional sets
collapse down to them. For example a decreasing sequence of spheres may have a
limit consisting of just a point. The n-1 dimensional space of outward
pointing vectors can be used for such a measure. It represents all $C^{2}%
$\ rectifiable sets of codimension at least one of $\mathbb{R}^{n}$ as
rectifiable sets in $\mathbb{R}^{n}$\textsf{X}$\mathbb{S}^{n-1}$ with n-1
dimensional Hausdorff measure. When viewed as $(n-1)$-rectifiable varifolds or
currents in $\mathbb{R}^{n}$\textsf{X}$\mathbb{S}^{n-1}$ they come equipped
with compactness theorems. The projection of their limits to $\mathbb{R}^{n}$
recovers the rectifiable sets of mixed dimensions giving the limits for the
desired topology on subsets of $\mathbb{R}^{n}$. Both varifold and current
compactness are required as there are sequences, such as honeycombs, that
converge as varifolds but not as currents. Conversely sequences such as lifts
of polyhedral approximations converge as currents but not as varifolds.

\end{abstract}
\maketitle
\tableofcontents

\section{Introduction}

This paper is the second of two, the first of which [MS] provides an example
(figure 1) of a physical optimization problem involving minimal surfaces and
threads of viscoelastic fluids with prescribed boundary and initial
conditions. This has a solution which is the limit of minimizers of an energy
functional, that is the limit of images of harmonic maps. The regularity
associated with minimizers, such as harmonic map images, can help establish
existence of limits of such minimizers. The limiting process in our example
has two features of interest. It involves dimension collapsing of the images
of the harmonic maps, where the lower dimensional sets are part of the
solution, and discontinuities can arise in bubbling. Although in our example
each limit of images can be achieved as a point-wise limit by careful
construction of a sequence of harmonic maps once it has been found [MS], a
general sequence will not have a defined point-wise limit.

This, the second paper provides two versions of a topology that overcomes both
the problem of requiring the lower dimensional sets to be kept in the limit,
and the problem of discontinuities due to bubbling in a sequence of harmonic
map images. This is done using sphere bundle measures (on $\mathbb{R}^{n}%
$\textsf{X}$\mathbb{S}^{n-1}$) to represent sets in $\mathbb{R}^{n}$ and to
take limits of sphere bundle measures rather than the underlying sets, then to
project down in the limit. The sphere bundle measure compactness has the added
advantage of providing a topology which compactifies certain unions of $C^{2}$
$j$-rectifiable subsets of $\mathbb{R}^{n}$ with $j$ varying from 0 to $n-1$,
where the boundaries and singular sets are also $C^{2}$ rectifiable.
Fortunately our motivating example consists of taking limits of sets which are
$C^{2}$ $(n-1)$-rectifiable, as they are images of harmonic maps. More
generally we would hope for the required regularity whenever we take sequences
of minimizers of functionals.

\subsection{Motivating example for this paper}

Figure 1 shows a sequence of images of harmonic maps from varying annuli into
$\mathbb{R}^{3}$ with fixed boundary circles [MS]. A topological cylinder can
be seen shrinking down to two discs union a straight segment. Such a sequence
could represent a physical optimization problem of a liquid film with
elasticity and surface tension forming a minimum energy configuration. The
surface tension causes area to be locally minimized resulting in the two
discs, while the elasticity causes the residual thread connecting the two
discs to remain.

\begin{center}%
{\includegraphics[
trim=0.008036in 0.056290in -0.008036in -0.056290in,
width=0.8867in
]%
{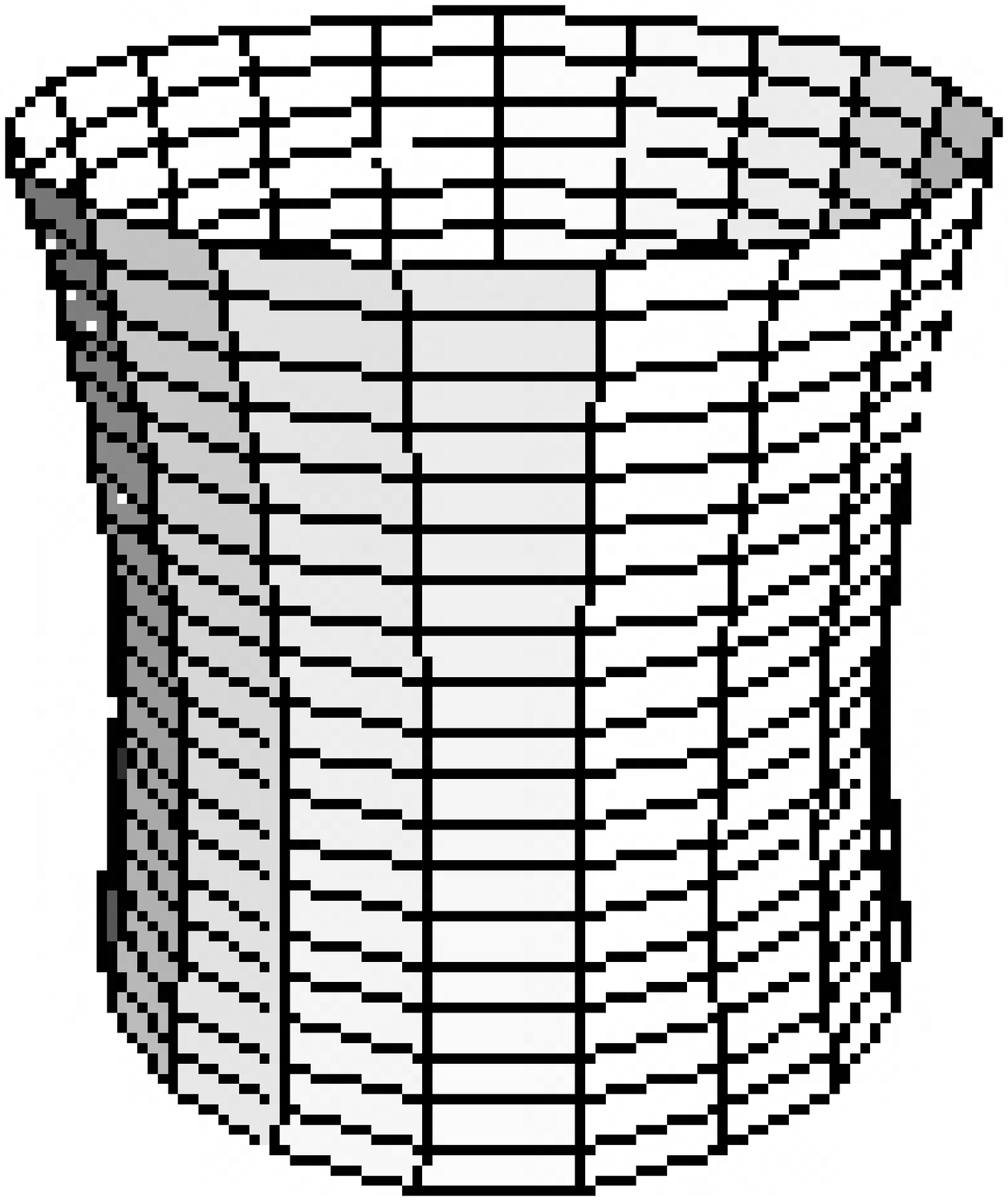}%
}%
{\includegraphics[
width=0.8867in
]%
{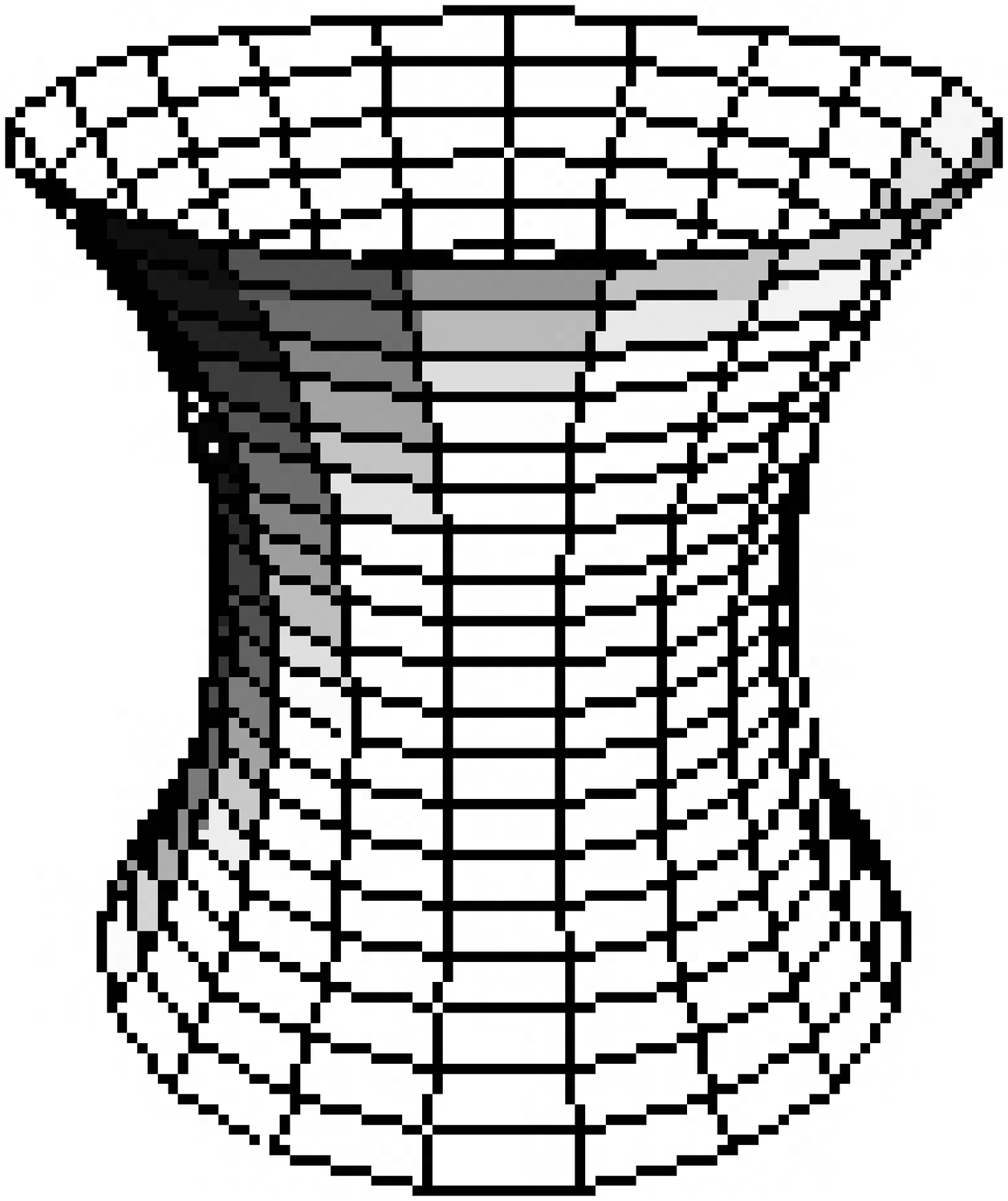}%
}%
{\includegraphics[
width=0.8867in
]%
{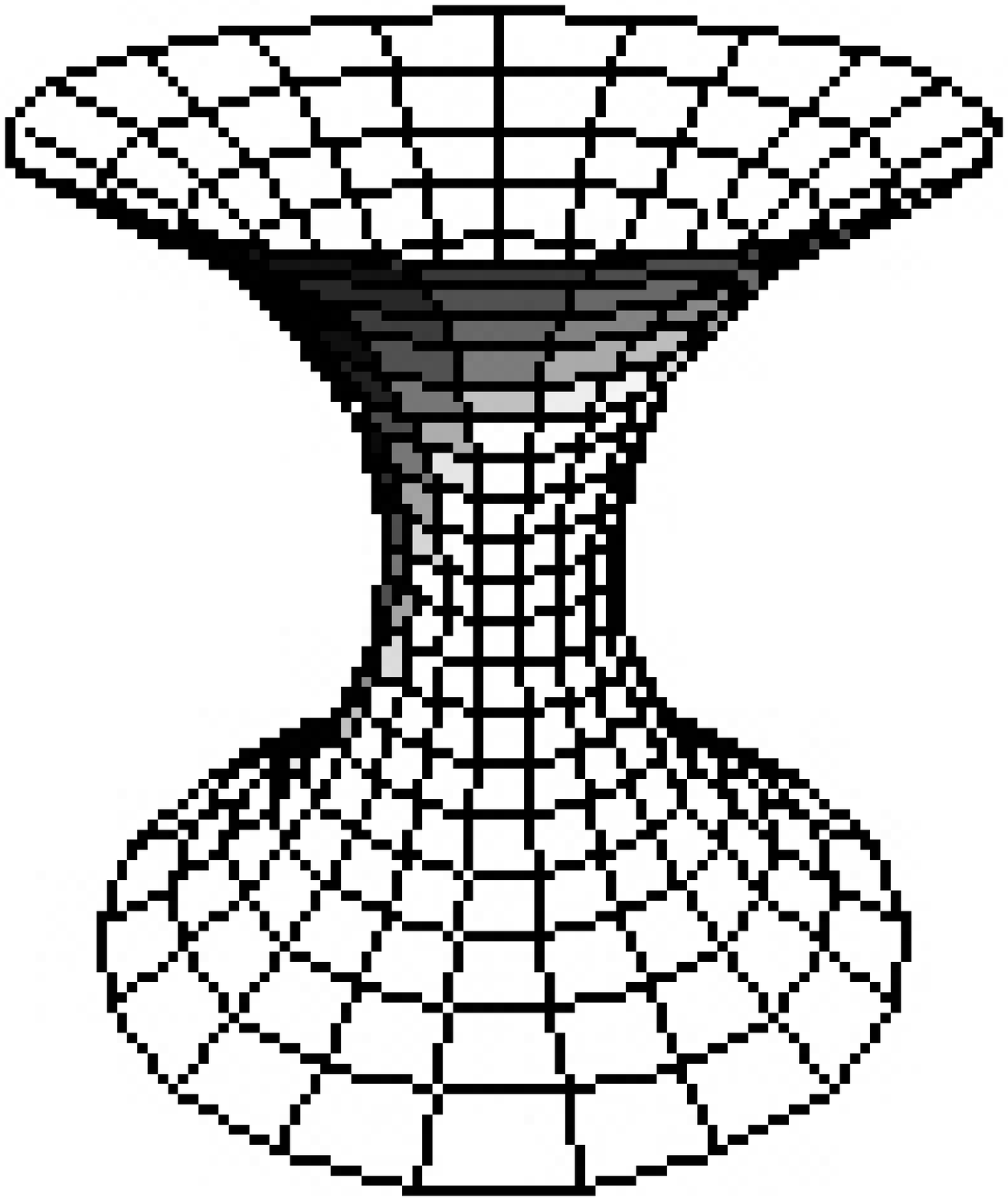}%
}%
{\includegraphics[
width=0.8867in
]%
{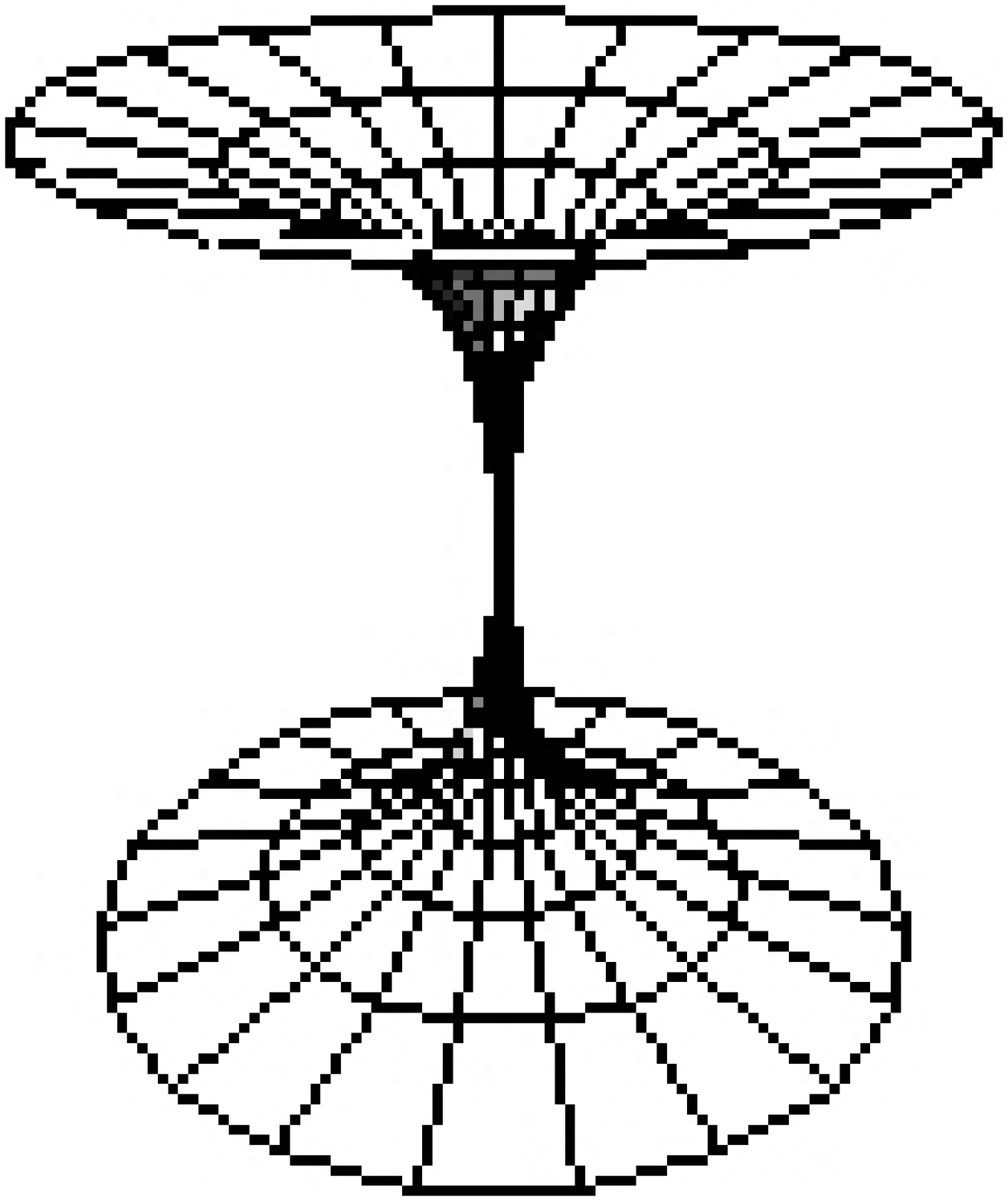}%
}%
{\includegraphics[
width=0.8925in
]%
{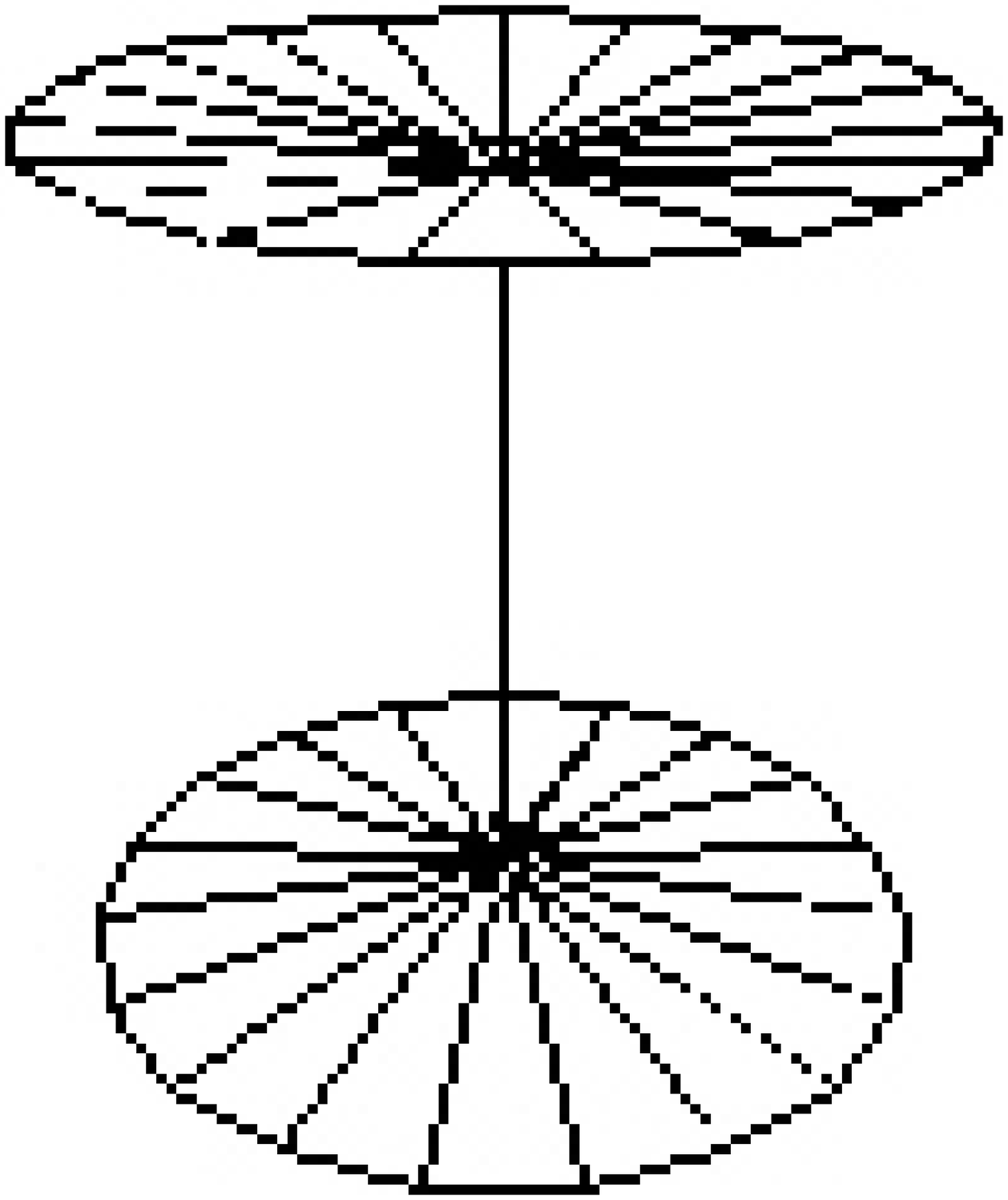}%
}%

\textbf{Figure 1: Dimension collapsing}
\end{center}

Conventional geometric measure theory provides measures and topologies with
compactness and regularity theorems that represent the above sequence in terms
of area measure. This will represent the cylinders and discs, but will miss
the residual straight line segment. Our goal is to find a measure and a
topology which will represent both the area measure and the straight line
segment in the limit. This would enable geometric measure theory compactness
and regularity theorems to be applied to these more general problems involving
sequences of sets whose dimensions collapse in the limit.

We take a measure based on the normal or outward pointing vectors to obtain
the desired concentration of measure around the straight line segments as
curvature increases and becomes infinite. In the limit we obtain a circle of
outward vectors concentrated at each point in the straight line segment, thus
giving 2-dimensional measure along the length the segment. We will see in
section 2 that outward vectors also allow boundary to be weighted, giving
another way for a measure to capture dimension collapsing, for example when an
interval shrinks in length down to a point. The two masses at each end of the
interval coincide in the limit, leaving mass at the point.

This new topology requires, among other uniform bounds conditions, the set of
outward vectors to be rectifiable. In our example application the sequence has
this regularity by virtue of being images of harmonic maps. So in general, in
applications we would take sequences of sets which are themselves minimizers
of some functional which endows them with the necessary regularity to have
rectifiable outwards vectors. Typically being $C^{2}$ almost everywhere, with
$C^{2}$ boundary or singular sets will be sufficient to ensure outward vector
regularity. See section 3 for necessary and sufficient conditions.

\subsection{Background}

Our approach and applications follow on from two approaches within the
mathematics of minimal surfaces. The use of harmonic maps to find minimal
surfaces [D] and the development of measure based compactness theorems, such
as for varifolds [Al] and for currents [FF], in geometric measure theory. For
introductions to the is field see [MF] for many examples, and [LY],[S] and
[HS], give deeper treatments, with [F1] as the general case reference not
written for the uninitiated. [Ma] covers sets from a measure theoretic
viewpoint. The set of outward pointing vectors on a smooth subset of
$\mathbb{R}^{n}$ of any dimension will be an n-1 manifold on the
$\mathbb{S}^{n-1}$ bundle of $\mathbb{R}^{n}$. Thus can be treated as an
$\left(  n-1\right)  $-rectifiable current or varifold in the $\mathbb{S}%
^{n-1}$ bundle, $\mathbb{R}^{n}$\textsf{X}$\mathbb{S}^{n-1}$. This enables
current or varifold compactness to be used on the $\mathbb{S}^{n-1}$ bundle so
that a limit is obtained which can be projected down to $\mathbb{R}^{n}$.

The underlying objects of geometric measure theory are rectifiable sets. These
have the following definitions:

\begin{definition}
Rectifiable sets

An m-dimensional set in $\mathbb{R}^{n}$ is countably n-rectifiable if

$M=M_{0}\underset{}{\overset{}{\cup}}\left(  \underset{j=1}{\overset{\infty
}{\cup}}F_{j}\left(  \mathbf{A}_{j}\subset\mathbb{R}^{m}\right)  \right)
,H^{m}\left(  M_{0}\right)  =0,$and\textit{\ }$F_{j}\left(  \mathbf{A}%
_{j}\right)  \rightarrow\mathbb{R}^{n}$, are Lipschitz functions for
j=1,2,3.... [S. p 58]
\end{definition}

This leads to the equivalent definition [S. p 59]

\begin{definition}
An n-dimensional set is countably n-rectifiable iff

$M\subset\underset{j=0}{\overset{\infty}{\cup}}N_{j\text{ }}$where
$H^{m}\left(  N_{0}\right)  =0$ and each $N_{j},$ $j\geq1$ is an m-dimensional
embedded $C^{1}$ submanifold of $\mathbb{R}^{n}$.
\end{definition}

Any rectifiable current or rectifiable varifold has underlying rectifiable set
on which the current or varifold has density
$>$
0. Varifolds should be considered as rectifiable sets which have a measurable
point-wise density [MF p11][F1, 2.9.12, 2.10.19], and unique approximate
tangent cones almost everywhere. See [MF p28] for approximate tangent cones
and [MF p92] for points without unique tangent cones. Currents, duals to
forms, should be considered as sets with density and unique approximate
tangent cones and orientation.

We also use the notion of $C^{2}$ $k$-rectifiable sets which were introduced
by Anzellotti [An].

\begin{definition}
$C^{2}$ $k$-rectifiable are sets for which the $N_{j}s$ in definition 1.2 are
embedded in $C^{2}$ submanifolds
\end{definition}

These $C^{2}$ submanifolds will have $C^{1}$ lifts, enabling the lifts of the
$N_{j}s$ to be rectifiable.

More specifically, our unions of rectifiable sets can be regarded as
multivarifolds as introduced by Dao and Fomenko[DF][Fo]. They in fact used
these multivarifolds to express limits of manifolds with dimension collapsing
just as in our example. [Fu][F2] show these can be the projections into
$\mathbb{R}^{n}$ of currents which can be seen as representing curvature
measure of the underlying sets. This is precisely the measure that
concentrates around parts of manifolds that collapse in dimension.

We can now state the main results of this paper:

\subsection{Statement of theorem 1.1}

This theorem gives the main gaol of the paper to provide a measure based
topology that captures dimension collapsing.

\begin{theorem}
Let M$_{i}$ be a sequence of unions of $C^{2}$ j-rectifiable sets of positive
H$^{j}$ measure in $\mathbb{R}^{n}$, (with j taking on different values;
$0\leq j\leq n-1$) where each of the outward vector lifts, $\widetilde{M}_{i}%
$, will all be n-1 rectifiable sets in $\mathbb{R}^{n}$\textsf{X}%
$\mathbb{S}^{n-1}$.

If the corresponding sequence $\widetilde{M}_{i}$ of outward vector lifts
converge as n-1 rectifiable currents or n-1 rectifiable varifolds in
$\mathbb{R}^{n}$\textsf{X}$\mathbb{S}^{n-1}$, then the image M of the
projection of the limit current or varifold, will itself be a union of
j-rectifiable sets in $\mathbb{R}^{n}$, ($0\leq j\leq n-1$). Also all
dimension collapsing will be captured. If the first variation of the lifts is
uniformly bounded then the lifts may converge as varifolds, if the homological
boundary is uniformly bounded then the lifts may converge as currents.
\end{theorem}

See sections 3 and 5 for conditions on when convergence as varifolds or
currents will occur.

\subsection{Outline of paper}

In section 2, we examine the process of taking the outward vectors for smooth
submanifolds with boundary of all dimensions in $\mathbb{R}^{n}$. These give
rise to a n-1 dimensional set in $\mathbb{R}^{n}$\textsf{X}$\mathbb{S}^{n-1}$,
the trivial sphere bundle over $\mathbb{R}^{n}$ . We see why we need to
represent more than just the vectors normal to tangent planes by including
extra vectors on boundary points. This ensures the measure captures dimension collapsing

Then in section 3 we examine these n-1 dimensional sets in $\mathbb{R}^{n}%
$\textsf{X}$\mathbb{S}^{n-1}$ as rectifiable n-1 varifolds, which gives them a
topology and compactness theorem. For compactness they need to meet strict
conditions, finite first variation, that would not be met for example by
successive polyhedral approximations to any smooth submanifold. We also
explore more fully the conditions on M such that the lift constructed in
section 2 will be rectifiable

Section 4 shows how sets in $\mathbb{R}^{n}$ can be obtained by projecting n-1
rectifiable varifolds from $\mathbb{R}^{n}$\textsf{X}$\mathbb{S}^{n-1}$ to
$\mathbb{R}^{n}$. This process is not the inverse of the process in section 2
as extra lower dimensional sets representing boundary or produced as a side
effect. Specifically, n-1 rectifiable varifolds project down to a union of
rectifiable varifolds of dimensions that can vary from 0 to n-1.

Section 5 shows why polyhedral approximations in $\mathbb{R}^{n}$ represented
as n-1 rectifiable varifolds in $\mathbb{R}^{n}$\textsf{X}$\mathbb{S}^{n-1}$
will not satisfy the conditions for varifold compactness, but as n-1
rectifiable currents they can satisfy conditions for rectifiable current
compactness. This is because the sets representing normal vectors to polyhedra
cannot be of finite first variation if polyhedral edge length is infinite.
However current compactness has no first variation condition, so with
appropriate `filling in' a polyhedron with infinite edge length but finite
total distributional curvature (edge length integral of dihedral angle) can be
represented as an n-1 rectifiable current.

However there is one drawback of current limits that varifolds do not have.
Mass can disappear in the limit due to cancellation of measures with opposite
orientations, or due to homothetic contraction of manifolds such as tori which
have degree zero Gauss maps. A method is given to overcome this using the
coarea formula to further fill in the lifts until each lift is a union of
$\mathbb{S}^{n-1}$s.

\section{The sphere bundle measure of submanifolds in $\mathbb{R}^{n}$}

We give a general construction (see figures 2 to 5) for lifts of a union of
submanifolds of mixed dimensions from 0 to $n-1$. These lifts will then be
represented as measures in $\mathbb{R}^{n}$\textsf{X}$\mathbb{S}^{n-1}$ where
geometric measure theory compactness theorems will give limits. We will not
address codimension zero subsets as these can be represented by taking
$\mathbb{R}^{n}$ as an isometrically embedded subset of $\mathbb{R}^{n+1}$ and
applying the construction within $\mathbb{R}^{n+1}$ instead.

\subsection{Definitions of outward vector lift}

\begin{definition}
An outward vector to a j-dimensional rectifiable set $M$ in $\mathbb{R}^{n}$
at a point $p$ is the set of unit vectors in the tangent space of
$\mathbb{R}^{n}$ which has a zero or negative inner product with any vector in
the approximate tangent cone of $M$ at $p$.
\end{definition}

In general ($H^{j}$ a.e.) this will be the normal vectors making up an $n-1-j
$ dimensional set in $\mathbb{R}^{n}$\textsf{X}$\mathbb{S}^{n-1}$. These unit
vectors can be represented as points in $\mathbb{S}^{n-1}$. Let $O(p)$ be the
set of all outward vectors at a point $p$, $O(p)=\{\mathbf{x}\in
\mathbb{S}^{n-1}:\mathbf{x}.\mathbf{v}\leq0,\forall\mathbf{v}\in X\}$ where
$X\subset\mathbb{S}^{n-1}$ is the link of the approximate tangent cone at $p$.

\begin{definition}
The \textbf{special case} lift of a set $M$ \textbf{with a unique approximate
tangent cone everywhere} is given by $\widetilde{M}=\{(p,O(p))\subset
\mathbb{R}^{n}$\textsf{X}$\mathbb{S}^{n-1}:p\in\overline{M}\}$ where
$\overline{M}$ is the closure of $M$. We can think of this as endowing open
sets as with the boundary of their closure.
\end{definition}

This definition only makes sense where $\overline{M}$ has a unique approximate
tangent cone [MF p28] [S] [F1]. Call $C$ the set of all points in
$\overline{M}$ with unique tangent cones. The approximate tangent cone of
$\overline{M}$ at $p$ can be thought of as a generalization of a linear
approximation to $\overline{M}$ at $p$. It allows for singularities such as
boundary or sections of planes coming together. For an m-submanifold the cone
will be the cone of a link which will be an m-1 manifold in the unit sphere in
the ambient space.

It is known that n-rectifiable sets have tangent spaces almost everywhere.
Therefore $\overline{M}-C$ is a set of zero measure, with respect to the
dimension of $\overline{M}$. See [MF p 92] for a point without a unique
approximate tangent cone.

So for$\ $all points, including those that do not have unique tangent cones we define.

\begin{definition}
$Oa(q)=\{(\mathbf{v}\in\mathbb{S}^{n-1}\mathbf{):v}\notin\{\mathbf{y}:\forall
r\exists x_{r}\in\overline{M}\cap B_{q,r}:\mathbf{y}.\overrightarrow{x_{r}%
-q}>0\}\}$
\end{definition}

\begin{definition}
The lift $\widetilde{M}=\{(p,O(p))\cup(q,Oa(q)):p\in C,q\in\overline
{M}-C\}\subset\mathbb{R}^{n}$\textsf{X}$\mathbb{S}^{n-1}.$
\end{definition}

Remark: In general $Oa(p)\neq O(p)$ where both exist. For example take a point
p on a circle in $\mathbb{R}^{2}$. $O(p)$ will consist of two normal vectors,
but $Oa(p)$ consists only of the outward pointing vector. In fact for a point
$p$ on the surface of a catenoid as in figure 1, $Oa(p)=\varnothing$.

Remark: In sets which are smooth embeddings of polygonal or simplicial
complexes, i.e.: the usual singular complexes, we do have a unique approximate
tangent cone everywhere.

Remark: Note that as the lifts are to be represented as varifolds, the closure
$Oa(p)$ is equivalent as a varifold to $Oa(p)$.

For technical reasons, to use current compactness, we also define a `filled
in' version of the lift, $\widetilde{M_{F}}$ to make the lift homologically
and topologically boundariless. Any remaining boundary on the lift is filled
in fiber by fiber. This means that there is a mass of maximum that of
hemisphere of the appropriate dimension added in each fiber. Also the lift is
filled in to ensure that the lift can be expressed as a union of spheres of
degree 1 under the projection map onto $\mathbb{S}^{n-1}$ for theorem 5.1.
This deals with cusps and non-manifold singularities where odd numbers of half
hyperplanes come together. Call these extra points $F$. Now we can define, for
use in thm 5.1 (and see figure 4):

\begin{definition}
The lift, $\widetilde{M}_{F}$= $\widetilde{M}\cup F$
\end{definition}

\begin{theorem}
: H$^{n-1}$($\widetilde{M}$)$\geq$H$^{n-1}$($\mathbb{S}^{n-1}$). For all
j-rectifiable M, j
$<$
n.
\end{theorem}

\begin{proof}
Take the convex hull of the $\overline{M}$. Every strictly convex point of the
hull will be in $\overline{M}$, because $\overline{M}$ is compact. Whether
such a point is in $C$ or $\overline{M}-C$, and will have an outward vector
lift locally containing the outward vector lift of the convex hull at those
points. Every unit vector in $\mathbb{R}^{n}$ (representable as a point in
$\mathbb{S}^{n-1}$) will occur in the outward vector lift of $p$ for at least
one $p$.
\end{proof}

\subsection{Examples of lifts}

For any compact boundariless smoothly embedded n-1 manifold $M$ in
$\mathbb{R}^{n}$.

$\widetilde{M}=\{(x,\mathbf{v}),(x,-\mathbf{v}),x\in M,\mathbf{v}\in
\mathbb{S}^{n-1},\mathbf{v\in}TM_{x}\bot\}$ $\subset$ $\mathbb{R}^{n}%
$\textsf{X}$\mathbb{S}^{n-1}$ See figure 2 for the lifts of a torus and a
sphere embedded in $\mathbb{R}^{3}$.%

\begin{center}
\includegraphics[
width=1.9933in
]%
{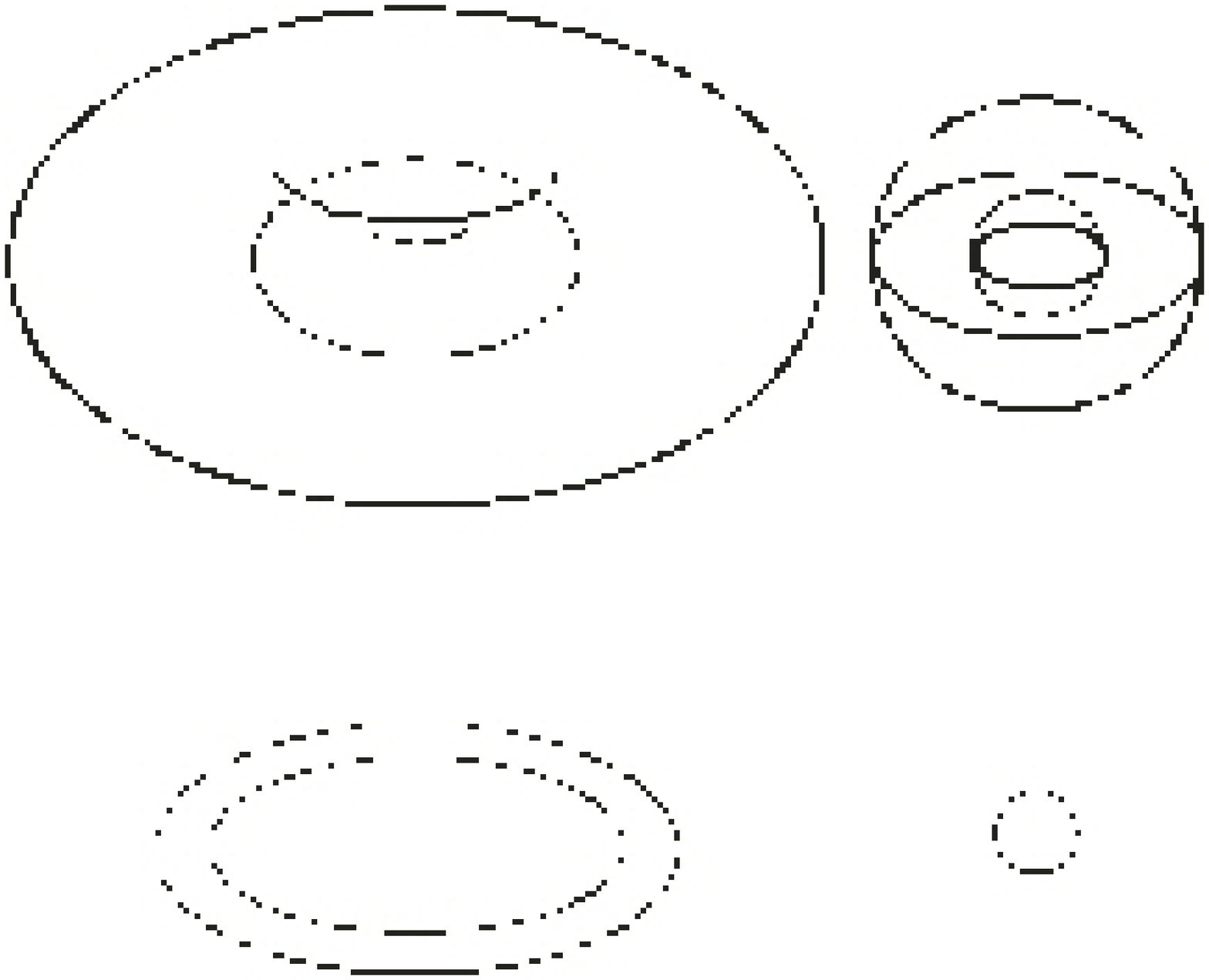}%
\end{center}

\begin{center}
\textbf{Figure 2: Lifts of normal vectors of }$\mathbf{M}$
\end{center}

Note that in each case $\widetilde{M}$ is a double cover of $M$.

We now show the outward vectors and the lift of an interval collapsing down to
a point to see how the boundary measure remains at the point.%

\begin{center}
\includegraphics[
width=3.1457in
]%
{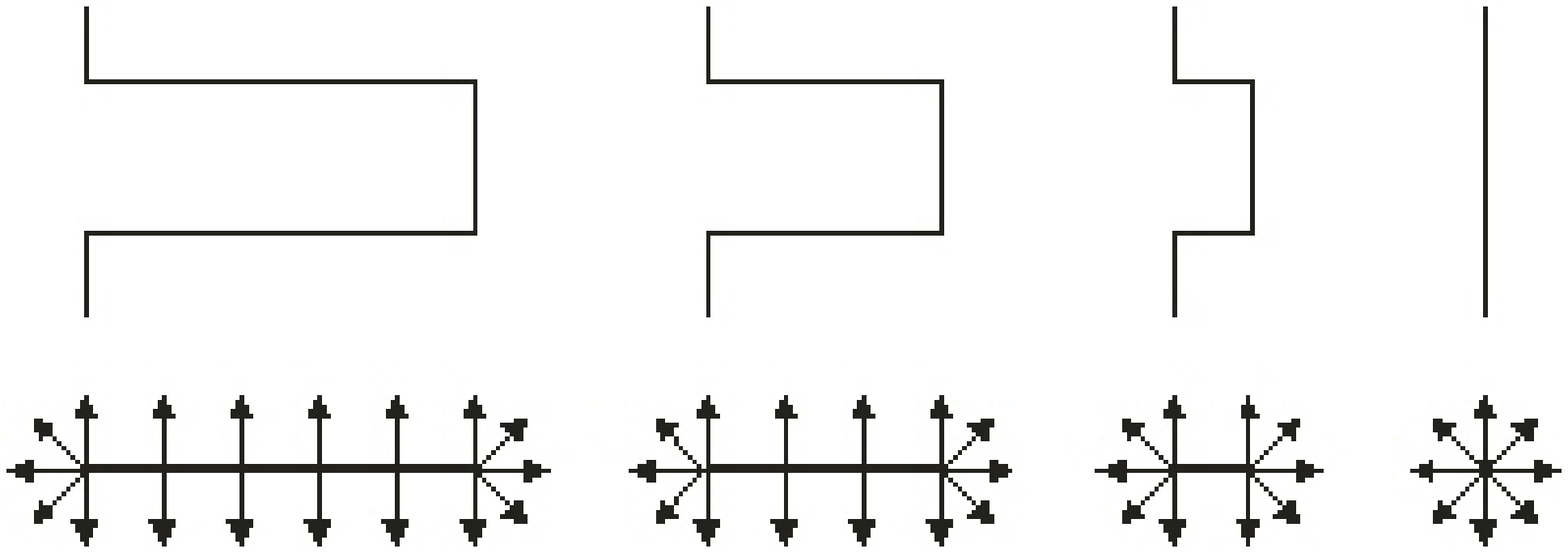}%
\end{center}

\begin{center}
\textbf{Figure 3: Lifts of outward vectors for dimension collapsing with
boundary.}
\end{center}

Figure 3 shows a sequence of line segments in $\mathbb{R}^{2}$ with their
outward pointing vectors, and the lifts of those vectors to $\mathbb{R}^{2}%
$\textsf{X}$\mathbb{S}^{1}$ are schematically represented above The
$\mathbb{S}^{1}$\ component is represented vertically.

In the case of a polygon collapsing, the outward pointing vectors at the
vertices would end up leaving measure at a point in the limit. See the outward
pointing vectors of the right angle as shown in figure 4. The vertical section
above the right angle would remain under homothetic contraction. This is in
addition to the vertical sections above the boundary points.%

\begin{center}
\includegraphics[
width=0.8501in
]%
{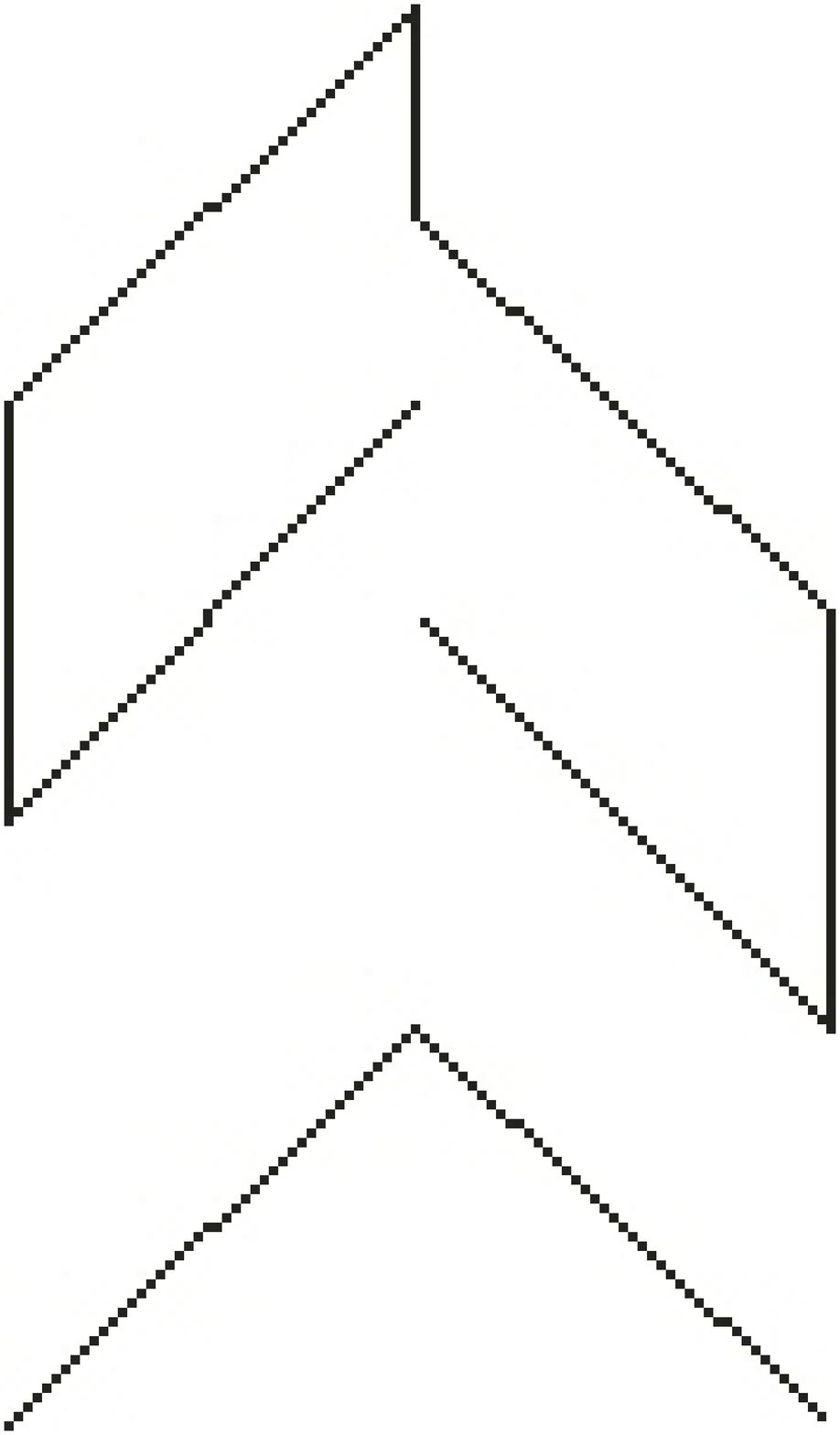}%
\end{center}

\begin{center}
\textbf{Figure 4: Lift of outward pointing vector of a right angle}
\end{center}%

\begin{center}
\includegraphics[
width=2.831in
]%
{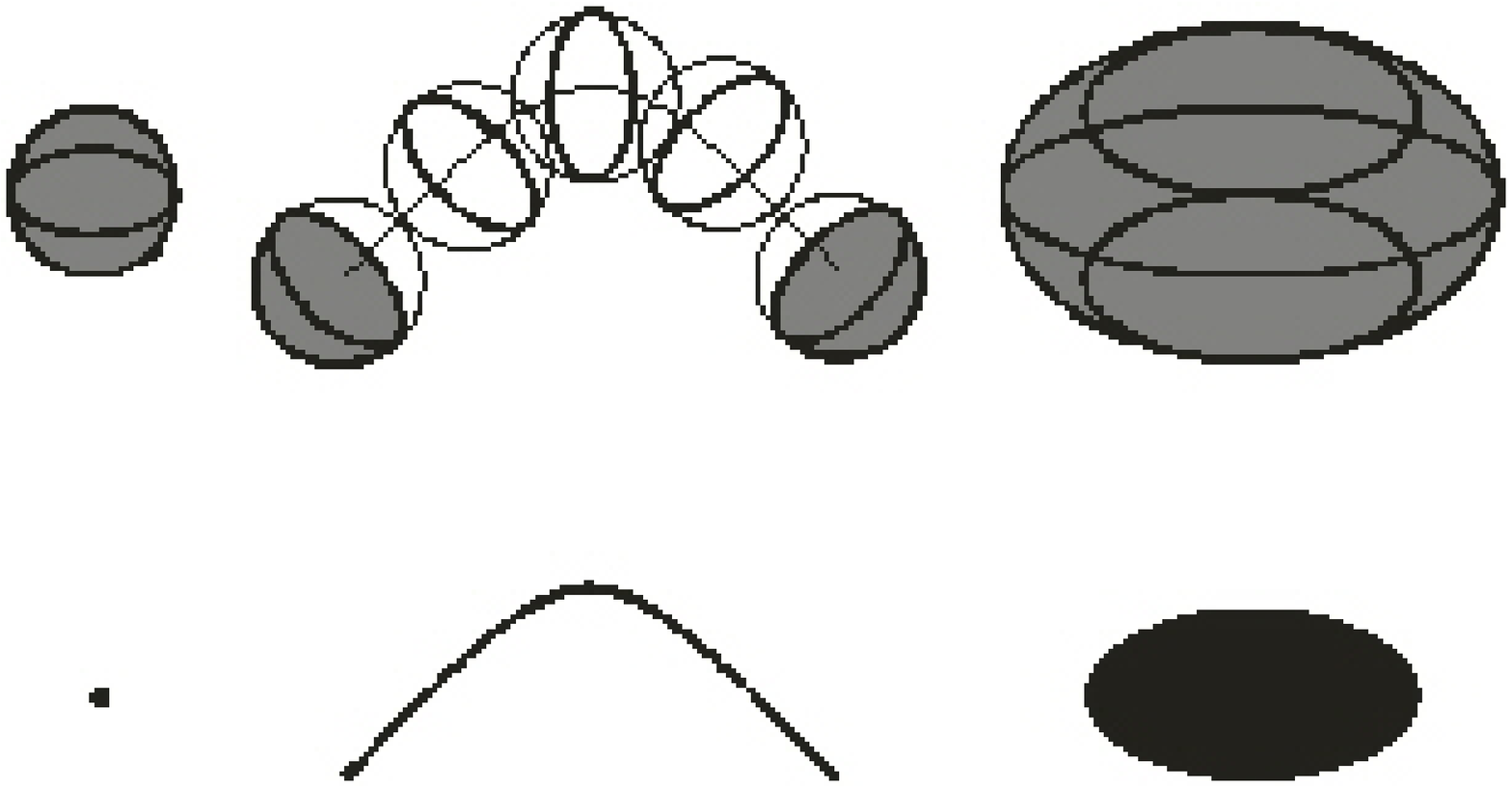}%
\end{center}

\begin{center}
\textbf{Figure 5: Outward vector lifts of a point a curve and a disc in
}$\mathbb{R}^{3}$
\end{center}

Figure 5 shows a different kind of schematic presentation of outward vector
lifts to $\mathbb{R}^{3}$\textsf{X}$\mathbb{S}^{2}$. Because the lift is in a
5-dimensional space, it has been smoothed out in the figure 5. The lifts are
not smooth as they appear, because in reality the lift of outward vectors at
the boundary meet the lift of vectors near the boundary orthogonally as in
figures 3 and 4. Topologically however, figure 5 gives the correct impression.

The point, the curve and the discs all have lifts that are topological
surfaces in $\mathbb{R}^{3}$\textsf{X}$\mathbb{S}^{2}$. As every vector is
outward from a point, the representation of the point is simply the whole
sphere. The points in the interior of the curve are represented by copies of
$\mathbb{S}^{1}$ as the normal vectors to the tangent space. Each boundary
point of the disc lifts to a semicircle of points.

\section{Varifold compactness for $\mathbb{S}^{n-1}$ bundle measures}

\subsection{Definition of varifold and total first variation}

The following is taken from Allard [Al].

A $k$-varifold on $N$, an smooth n-manifold, is a radon measure on the bundle
over $N$ whose fiber at each point $p$ of $N$ is the Grassmann manifold of $k$
dimensional hyperplanes in the tangent space to $N$ at $p$. A rectifiable
integral $k$-varifold is strongly approximated by a positive integral linear
combination of varifolds corresponding to continuously differentiable
$k$-dimensional submanifolds of $N$.

$IV_{k}(U)$ is the set of $k$-rectifiable integral varifolds on $U,$a open
subset of N. $|V|$ is $k$-dimensional Hausdorff measure of the $k$-dimensional
set $M\subset U\subset N$ corresponding to $V$. $||\partial V||$ is the total
first variation of the $k$-dimensional Hausdorff measure of $M$.

Note: $\Vert\partial V\Vert=%
{\displaystyle\int\limits_{\widetilde{M}_{V}}}
\left\Vert \overrightarrow{H}\right\Vert dH^{k}+%
{\displaystyle\int\limits_{\partial\widetilde{M}_{V}}}
dH^{k-1}$, is the total first variation where $\overrightarrow{H}$ is the mean
curvature vector of $M$, in the distributional sense in the integral as
described in the proof.

\subsection{Varifold compactness}

\begin{theorem}
Suppose $U$ is an open subset of $\mathbb{R}^{n}$, $G_{1}$, $G_{2}$, ... are
open subsets of $U$, $U=\overset{\infty}{\underset{i=1}{\cup}}G_{i}$ and
$l_{1}$, $l_{2}$, ... are nonnegative real numbers. Then

$IV_{k}(U)\cup\{V:(\Vert V\Vert+\Vert\partial V\Vert)(G_{i})\leq l_{i}\}$ is compact.
\end{theorem}

This means that the space of uniformly finite mass and finite first variation
$k$-rectifiable integral varifolds is compact. Intergal varifolds have
positive integer density almost everywhere. So we now need to determine which
sequences in $\mathbb{R}^{n}$ will have lifts that converge as varifolds.

\subsection{Conditions on sets in $\mathbb{R}^{n}$ so that the lifts converge
as varifolds}

Conditions (3.3.i) and (3.3.ii) also apply for current compactness, but we
will see in section 5 that condition (3.3.iii) changes for currents.

\textit{(3.3.i) Sets and their boundaries, and lower dimensional faces, are
}$C^{2}$\textit{-rectifiable}. This will ensure that the lifts are rectifiable
as long as they have finite mass.

\textit{(3.3.ii) Sets and their boundaries, and lower dimensional faces, have
uniformly bounded finite mass and finite integrals of principle curvatures
(including distributional).} This keeps the mass of the lifts uniformly
bounded. Notice that we need to use distributional curvature. Consider the
example of a polygon in $\mathbb{R}^{3}$, each point on the face lifts to two
points in $\mathbb{R}^{3}$\textsf{X}$\mathbb{S}^{2},$ each point on an edge
contributes a semicircle and each vertex contributes a portion of a sphere
depending on the exterior angle of the polygon.

In general for a set $M$ of dimension $j$ the curvature bounds affect the
Jacobian [S: section 8] of the lift.%

\begin{equation}
\underset{1\leq k\leq j}{\Pi}\sqrt{1+\rho_{k}^{2}}%
\end{equation}

where $\rho_{k}$ is a principle curvature in a direction $k$ making up a local
basis of $TM,$ the tangent space of $M$ in $\mathbb{R}^{n}$.

Consider the lift to the sphere bundle of each of the $\mathit{j}$ directions
in $TM$. For each $k$ the lift has an identity component in an $\mathbb{R}%
^{n}$ direction, and a component in the $\mathbb{S}^{n-1}$ direction. Each
directional derivative therefore has a component of 1 and a component equal to
the curvature of $M$ in the $k$ direction. The lift of the outward vectors is
in fact the space of orthogonal vectors to the lift in the fiber
$\mathbb{S}^{n-1}$.

For the interiors of the manifolds and the boundaries this gives;

\begin{center}%
\begin{equation}
\underset{j}{\sum}\left(  \underset{M_{j}}{\int}\varpi^{n-j-1}\underset{1\leq
k\leq j}{\Pi}\sqrt{1+\rho_{k}^{2}}dH^{j}+\underset{\partial M_{j}}{\int}%
\frac{\varpi^{n-j}}{2}\underset{1\leq k\leq j-1}{\Pi}\sqrt{1+\rho_{k}^{2}%
}dH_{{}}^{j-1}\right)  <L
\end{equation}

\end{center}

where $M_{j}$ is the union of j-dimensional manifolds of $M$, and $\varpi^{n}$
is the volume of the unit $n$ sphere. So for the sequence $M_{i}$ this must be
a uniform bound. Now there is an additional term for interior non-smooth
points. This can be given by

\begin{center}%
\begin{equation}
\underset{j}{\sum}\underset{}{\underset{S\subset M_{j}}{\int}\varpi^{n-j}%
}KdH^{j-1}<L
\end{equation}

\end{center}

where $S$ is the non-smooth singular set, and $K$ represents the
distributional curvature on that set. In the case of $j=2$, $K$ is just the
exterior dihedral angle.

\textit{(3.3.iii) Sets have bounds on the integrals of the derivatives of
curvature.}

The lifts must have uniformly finite first variation. This requires bounds on
the integrals of the derivatives of curvature on the $M_{i}$ and finite
boundary mass. Also there must be finite mass on singular sets such as
polyhedral dihedral angle. This is because, as illustrated in figures 3 and 4
the lifts of non-smooth points will be non-smooth. Furthermore the lift of a
small dihedral angle in the case of $j=2$, will involve dihedral angles in the
lift equal to $\frac{\pi}{2}$. In other words lifting makes first variation
worse, by an unbounded factor. See section 5 on currents to see how to combat
this drawback of varifolds compactness.

\section{Projecting sphere bundle measure limits to sets in $\mathbb{R}^{n}$}

The limiting process in the previous section gives rise to limit with an
$(n-1)$-dimensional Hausdorff measure in $\mathbb{R}^{n}$\textsf{X}%
$\mathbb{S}^{n-1}$. When we project down to $\mathbb{R}^{n}$ locally we may
obtain anything from zero to $n-1$ dimensional Hausdorff measure depending on
the rank of the projection map. The rank of the map at a point depends upon
how many basis vectors of the tangent space as subset of $\mathbb{R}^{n}%
$\textsf{X}$\mathbb{S}^{n-1}$ lie in the $\mathbb{R}^{n}$ directions and how
many lie in the $\mathbb{S}^{n-1}$. In fact the rank is just the number of
basis vectors that lie in the $\mathbb{R}^{n}$ directions.

Accordingly say we have a limit varifold , we can decompose it into components
each of which is $n-1$ rectifiable and has a fixed rank under projection to
$\mathbb{R}^{n}$.

\begin{theorem}
An (n-1)-rectifiable set $A$ in $\mathbb{R}^{n}$\textsf{X}$\mathbb{S}^{n-1}$
can be decomposed into n-1 different (n--1)-rectifiable sets, $A_{i}$, ($0\leq
i\leq n-1),$ where the tangent space of $A_{i}$ at each point intersects the
tangent space of $\mathbb{S}^{n-1}$ in i-dimensional space.
\end{theorem}

\begin{proof}
This can be shown using the fact that an $(n-1)$-rectifiable set can be
defined in terms of being a subset of a countable union of $C^{1}$
$(n-1)$-dimensional submanifolds (Def. 1.2). It is clear that each subset of
fixed rank under the projection map is itself $(n-1)$-rectifiable as it is
contained in the original union of submanifolds.
\end{proof}

Now each of these sets has a projection so we need to show that each
projection is rectifiable of the correct dimension. We will prove this by contradiction.

\begin{theorem}
Each of the $A_{i}$s above projects down to an $(n-1-i)$-rectifiable set
$p(A_{i})\in\mathbb{R}^{n}.$
\end{theorem}

\begin{proof}
First we show $p(A_{i})\in\mathbb{R}^{n}$ is rectifiable as it is a projection
of a rectifiable set. As the projection of a $C^{1}$ submanifold is a $C^{1}$
submanifold then the projection is a rectifiable set using definition 1.2.
\end{proof}

The coarea formula [F1: 3.2.32, 3.2.31] relates the measures of the varifolds
and their projections up to sets of measure zero in the lift and in the
projection. Therefore we can now proceed to write down varifolds in
$\mathbb{R}^{n}$ in terms of the limit of lifts.%

\begin{center}
\includegraphics[
width=3.4014in
]%
{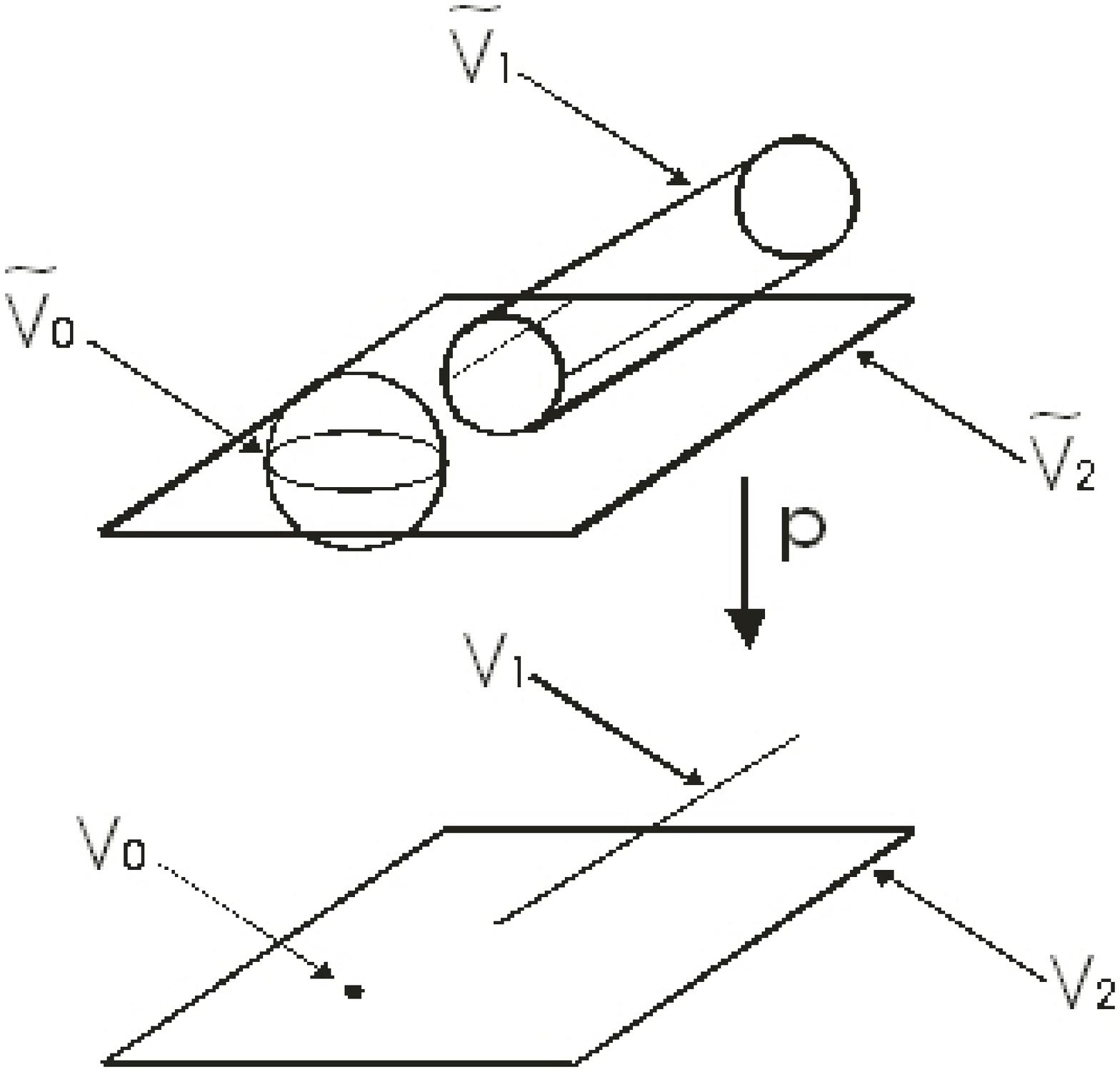}%
\end{center}

\begin{center}
\textbf{Figure 6: Decomposition of part of a limit varifold and projections}
\end{center}

Figure 6 shows how a part of a limit of lift varifolds $\widetilde{V}$ of
density 1 can be represented as a varifold:

\begin{center}
$\widetilde{V}=\overset{n-1}{\underset{i=0}{\sum}}\widetilde{V}_{i}$ in the
case $n=3$.
\end{center}

Also the projections are shown in $\mathbb{R}^{3}$ of the underlying sets of
the lift varifolds, which become the underlying sets for $V_{0}$, $V_{1}$,
$V_{2}$. As the point underlying $V_{0}$ and part of the line segment
underlying $V_{1}$ are contained in the planar region underlying $V_{2}$, we
can see that the limit of lifts in $\mathbb{R}^{n}$\textsf{X}$\mathbb{S}%
^{n-1}$ contains more information than just the set.

We can define the $V_{i}$s in alternate ways. First consider the following
simple equation that defines varifolds $V_{i}$ in terms of varifold
$\widetilde{V_{i}}$.%

\begin{equation}
V_{i}\left(  f\right)  =\frac{1}{\varpi^{n-1-i}}\underset{x\in\widetilde
{V}_{i}}{\int}f\circ pdH^{n-1}\left(  \mathbb{R}^{n}\text{\textsf{X}%
}\mathbb{S}^{n-1}\right)
\end{equation}

Recall $\varpi^{n}$ is the volume of the unit n sphere.

We shall now give a construction for the $Vi$s only terms of $\widetilde{V}$.
Using $\theta_{\widetilde{V}}$ to represent the density of $\widetilde{V}, $
which will be integer multiplicity almost everywhere, we can use modified
versions of 4.1 to obtain:%

\begin{equation}
V_{0}\left(  f\right)  =\frac{1}{\varpi^{n-1}}\underset{x\in\mathbb{R}^{n}%
}{\int}f(x)\left(  \underset{\left(  p^{-1}\left(  x\right)  \cap\widetilde
{V}\right)  }{\int}\theta_{\widetilde{V}}(p^{-1}\left(  x\right)
)dH^{n-1}(\mathbb{S}^{n-1})\right)  dH^{0}\left(  \mathbb{R}^{n}\right)
\end{equation}
\qquad%
\begin{equation}
V_{1}\left(  f\right)  =\frac{1}{\varpi^{n-2}}\underset{x\in\mathbb{R}^{n}%
}{\int}f(x)\left(  \underset{\left(  p^{-1}\left(  x\right)  \cap\widetilde
{V}\right)  }{\int}\theta_{\widetilde{V}}(p^{-1}\left(  x\right)
)dH^{n-2}(\mathbb{S}^{n-1})\right)  dH^{1}\left(  \mathbb{R}^{n}\right)
\end{equation}
\qquad%
\begin{equation}
V_{i}\left(  f\right)  =\frac{1}{\varpi^{n-1-i}}\underset{x\in\mathbb{R}^{n}%
}{\int}f(x)\left(  \underset{\left(  p^{-1}\left(  x\right)  \cap\widetilde
{V}\right)  }{\int}\theta_{\widetilde{V}}(p^{-1}\left(  x\right)
)dH^{n-1-i}(\mathbb{S}^{n-1})\right)  dH^{i}\left(  \mathbb{R}^{n}\right)
\end{equation}

We can show why each double integral picks out just the measure that is
desired. In finding $V_{j}$, the inner integral takes $H^{n-1-j}$ measure.
This will eliminate any contribution to a $V_{k}$ for $k>j$, as the fibers
above each point will have zero $H^{n-1-j}$ measure. Similarly the outer
integral will eliminate contributions to a $V_{k}$ for $k<j$. This is because
although on points in $supp(V_{k})$ a fiber will have infinite $H^{n-1-j}$
measure, this will only occur on a set zero $H^{j}$ measure on $\mathbb{R}%
^{n}$

Notice that for varifolds of dimension less than $n-1$ their underlying
rectifiable set can be infinite. This is because there is no lower bound on
the volumes of mass in each fiber. This is caused for example with a
polyhedral approximation. The lengths of edges can go to infinity. So these
are, in a sense, false projections as they correspond not to dimension
collapsing such as threads in figure 1 but curvatures. Now [F1 3.2.31] allows
for a cut off to be made so that only a subset of the projection is taken,
that is those points that have more than a critical lower bound of mass. This
new projection is then still rectifiable, and will then have finite mass. Call
these varifolds, $V_{i}\mid_{m},$.that represent varifolds whose lift fibers
all have at least mass $m$.

We have now proven part of theorem 1. The case where the first variation of
the lifts of a sequence of sets is uniformly bounded. In the next section we
will see that this is a limiting and avoidable restriction on sets, as long as
certain types of non-manifold point singularities are bounded.

Remarks:

\textbullet\qquad According to the area and coarea formulae there will be a
discrepancy between 4.1 on the one hand and 4.2 to 4.4 on the other. The map
$p^{-1}$ always has a Jacobian greater than one (by virtue that $p$ is a
projection) and less than infinity (because the measure of $\widetilde{V}$ is
finite). As we are not using density information on the $V_{j}$s, this
discrepancy is not important.

\textbullet\qquad The final stage is to take the support of each $V_{j}$ to
obtain the limit rectifiable set as desired. We may not get integer
multiplicity of the varifolds $V_{j}$, so the underlying rectifiable set is a
better canonical limit object, rather than the measure. One set can represent
everything, but as measures, $n-1$ separate measures are needed to represent
all of the set in $\mathbb{R}^{n}$ , although in the sphere bundle one measure
can represent the union. In particular varifolds representing boundary will
have density of
${\frac12}$
in general when no dimension collapsing occurs.

\textbullet\qquad In general the geometric relationship between the sets in
$\mathbb{R}^{n}$ and their lifts to outward pointing vectors in $\mathbb{R}%
^{n}$\textsf{X}$\mathbb{S}^{n-1}$ may break down in the limit. Under the right
conditions in the sequence of sets in $\mathbb{R}^{n}$ , such as smoothly
embedded connected submanifolds without boundary and with uniformly bounded
curvature everywhere, we might expect the limit of the lifts in $\mathbb{R}%
^{n}$\textsf{X}$\mathbb{S}^{n-1}$ to represent the outward pointing vectors of
its projection into $\mathbb{R}^{n}$ . This would enable us to carry over
regularity information from the sequence of sets in $\mathbb{R}^{n}$ to its
limit under our process. We can conjecture that such theorems can be proven
using basic geometric measure theory on the tangent space of the lift
varifolds, but are beyond the scope of this paper.

\section{Currents}

\subsection{The Polyhedral Case}

Here we explore the consequences of the need for a uniform bound on the
integrals of derivatives of curvature in order for the first variation of the
lifts to be uniformly bounded for varifold compactness. The worst case for
derivatives of curvature are where the curvature concentrates around small
regions rather than being evenly distributed over a manifold. In the limit
this gives polyhedra, so we will now explore the lifts of polyhedra and see
how to deal with them as measures in the sphere bundle.

\begin{center}%
{\includegraphics[
trim=0.000000in 0.000000in -0.001059in 0.000000in,
width=3.5566in
]%
{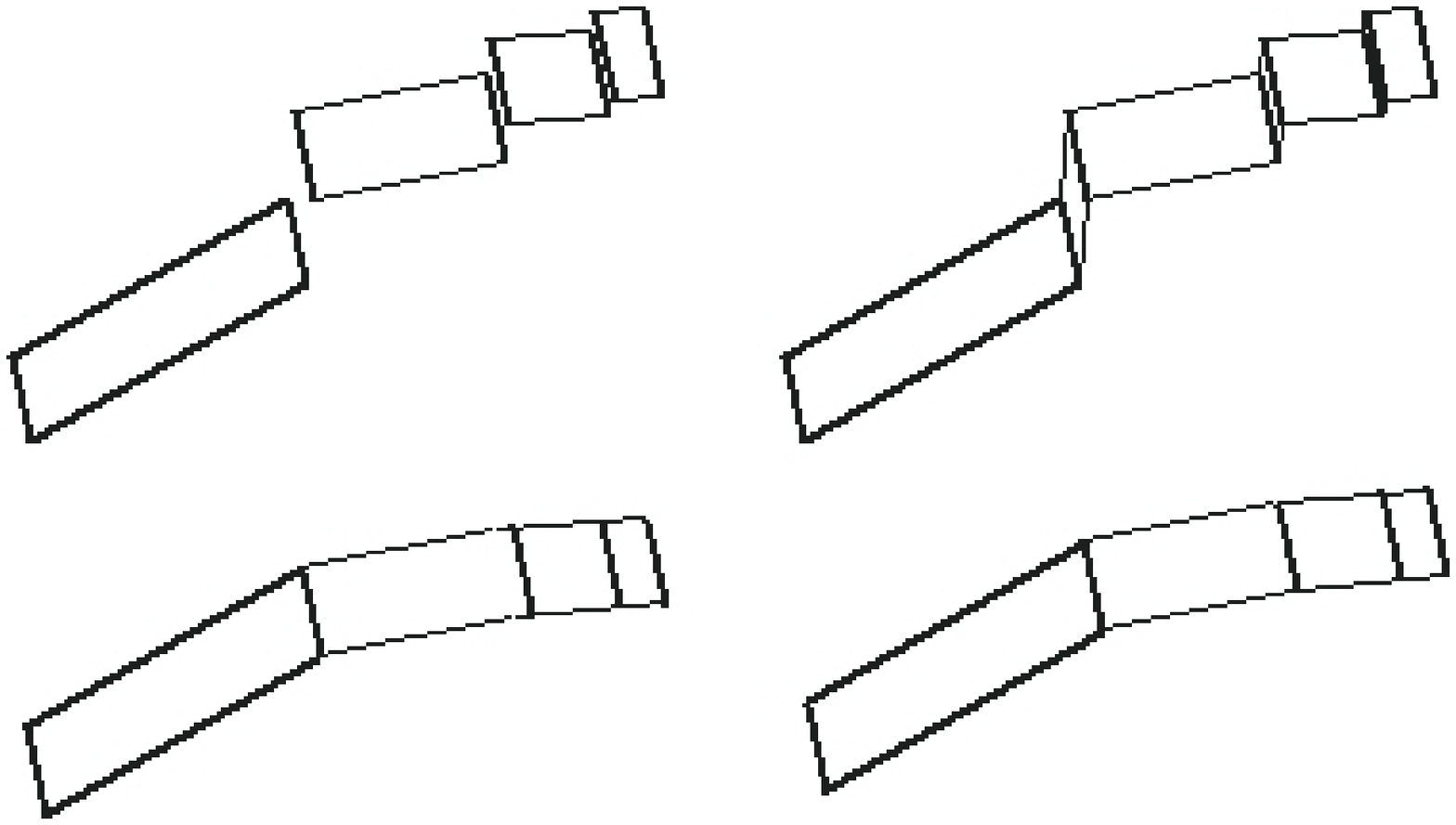}%
}%
\qquad%
{\includegraphics[
width=1.0552in
]%
{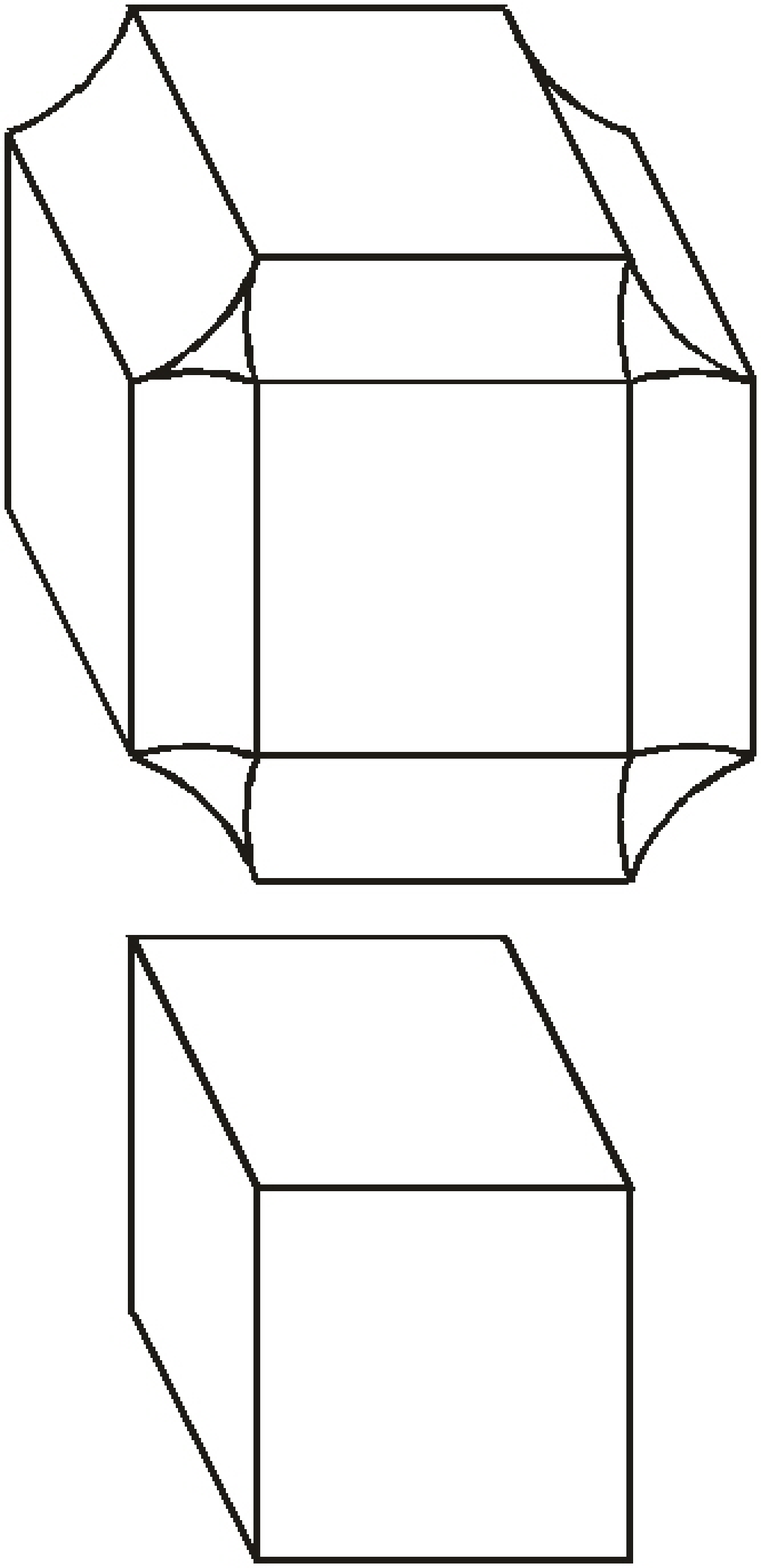}%
}%

\textbf{Figure 7: Polyhedra and their lifts}
\end{center}

Figure 7 shows two copies of a polyhedral strip and a cube in $\mathbb{R}^{3}
$ below with graphical representations of their various lifts above. Top left
shows a normal vector lift of the polyhedral strip, showing that two boundary
components lie in the fiber above an edge. Top middle shows how outward
vectors at each edge lift to fill in the discontinuities in the lift on the
left removing much of the new boundary. Note here that the fiber above every
point in the edge is filled in between the directions of the adjacent faces.
In fact the filled in part meets lifts of the adjacent edges orthogonally.
This adds first variation to the lift proportional to the length of the edges.
This means that varifold compactness in section 3 will not work for sequences
polyhedral approximations of smooth surfaces whose edge lengths will tend to infinity.

\subsection{Definition of Rectifiable currents}

An $n-$dimensional rectifiable current in $U$ is a continuous linear
functional on $D^{n}(U)$, the space of smooth n-forms with compact support in
$U$. These currents are denoted $D_{n}(U)$. For $n>0$ the currents can be
viewed as generalizations of the n-dimensional oriented submanifolds $M$
having locally finite $H^{n}$-measure in $U$. Given such an $M$ with an
orientation $\xi$, we can write $\xi(x)=\pm\tau_{1}\wedge...\wedge\tau
_{n}.\forall x\in M$ where $\tau_{1}\wedge...\wedge\tau_{n}$ is an orthogonal
basis for $T_{x}M$. Additionally, just as with rectifiable $n-$varifolds we
can endow $M$ with a density function $\theta(x)$, a positive locally
integrable function on $M$. Then we have the corresponding n-dimensional
current $T$ where:

$T\left(  \varpi\right)  =%
{\displaystyle\int\limits_{M}}
\left\langle \varpi\left(  x\right)  ,\xi(x)\right\rangle \theta(x)dH^{n},$
$\varpi\in D^{n}\left(  U\right)  $

We can define a boundary current $\partial T(\varpi)=T\left(  d\varpi\right)
$ in accordance with Stokes's theorem. When $\theta(x)$ is integer valued
$H^{n}$ almost everywhere we saw that the current is integer multiplicity.

\subsection{Statement of current compactness theorem}

The Federer-Fleming current compactness theorem [S],[FF] is as follows:

\begin{theorem}
Suppose $T_{j}$ is a sequence of integer multiplicity rectifiable n-currents
in $\mathbb{R}^{n+k}$ and $\partial T_{j}$ are integer multiplicity
rectifiable n-1 currents in $\mathbb{R}^{n+k}$ and $W$ is a compact subset of
$\mathbb{R}^{n+k},$ $sup_{j\geq1}(M_{W}(T_{j})+M_{W}(\partial T_{j}))<\infty$
then there will be a subsequence that converges weakly in $W $ to a
rectifiable current $T$, which is also an integer multiplicity rectifiable
$n$-dimensional current. $M_{W}\left(  T\right)  $ is the mass of the current
$T$ in a compact set $W$ and is given by $M_{W}(T)=\underset{\varpi}{\sup
}\left\vert T\left(  \varpi\right)  \right\vert ,\varpi\in D^{n}\left(
W\right)  ,\left\vert w(x)\right\vert \leq1,\forall x\in W.$ Weak convergence
is in the obvious sense, $T_{j}\rightarrow T\Leftrightarrow\underset
{}{\underset{j\rightarrow\infty}{\lim}\left(  T_{j}\left(  \varpi\right)
\right)  =T\left(  \varpi\right)  ,}\forall\varpi\in D^{n}$.
\end{theorem}

\subsection{Conditions on sets in $\mathbb{R}^{n}$ so that the lifts converge
as currents}

Current compactness differs from varifold compactness in one key way. First
variation is not required to be bounded, but homological boundary volume is.
Homological boundary can exist locally with no first variation, for example
see appendix. This occurs when orientation is assigned and for example half
spaces come together, such as a Y-singularity. The induced boundaries on each
of the three parts where they cannot all cancel.

The filled in lifts do satisfy current compactness conditions if the unfilled
lifts do. The filling in contributes mass according to the distributional
curvature of the $M_{i}$, in the surface case. So the hypotheses of bounded
curvature integrals including distributional curvature will suffice to ensure
current compactness.

We need conditions 3.3.i and 3.3.ii, but instead of 3.3.iii we need to add a
hypothesis that the topologically singular set, such as Y-singularities must
have a uniformly finite mass lift. However it weakens a hypothesis. The
geometrically singular set (approximate tangent cones have a dihedral angle)
can now be infinite as long as the integral of dihedral angle is uniformly
bounded as is the case with certain sequences of polyhedral approximation(this
was covered by 3.3.ii). Now we can write conditions on sets for the lifts to
converge as currents:

\textit{(5.4.i) Sets and their boundaries, and lower dimensional faces, are
}$C^{2}$ \textit{rectifiable}.

\textit{(5.4.ii) Sets and their boundaries, and lower dimensional faces, have
uniformly bounded finite mass and finite integrals of principle
curvatures(including distributional).}

\textit{(5.4.iii) The topologically singular set, such as Y-singularities,
must have a uniformly finite mass lift.}

Note that for unions of non-singular manifolds current compactness is easier
to achieve in the lifts than varifold compactness.

\subsection{The problem of mass cancellation in current compactness}

Say we have a sequence of oriented codimension 1 manifolds, with hypotheses
fulfilled. We can see two simple ways mass can be cancelled in the limit. One
is by two submanifolds converging to each other on a set of positive measure
with opposite cancelling orientations. Another is under homothetic
contraction. Take a torus in $\mathbb{R}^{3}$ with an orientation. Allow the
lift to have the same orientation. Under homothety the limit of the lift will
correspond to the push forward of the current under the Gauss map, as all
vector directions are represented in the fiber above the center of the
homothety in $\mathbb{R}^{3}$. As the degree is zero this will mean that the
current will cancel to zero in the limit. This does not give us any record of
dimension collapsing, in terms of residual measure. The next section gives a
method to overcome this problem under certain circumstances.

\subsection{Elimination of mass cancellation}

\begin{center}%
{\includegraphics[
width=1.4861in
]%
{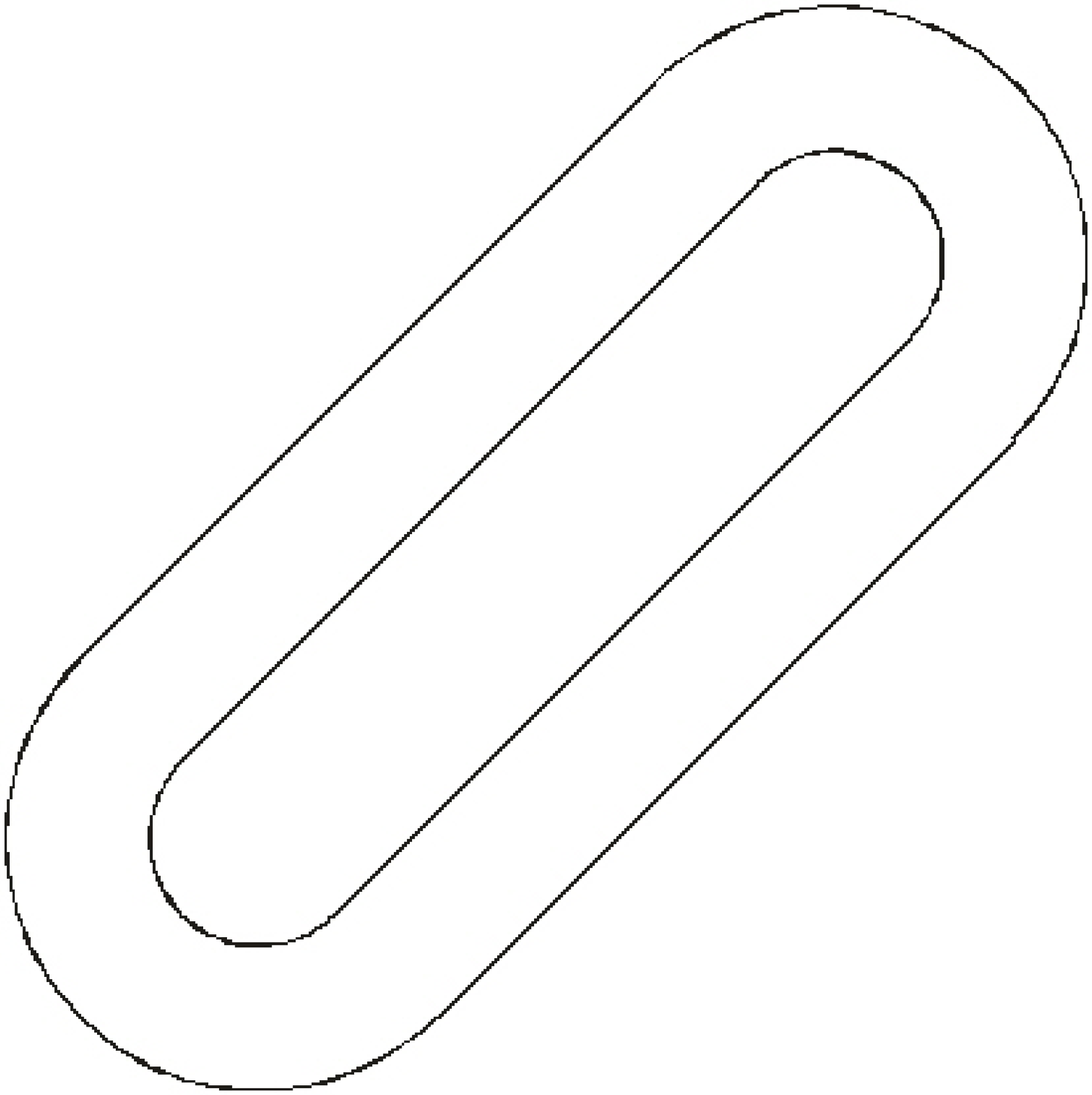}%
}%
{\includegraphics[
width=1.4861in
]%
{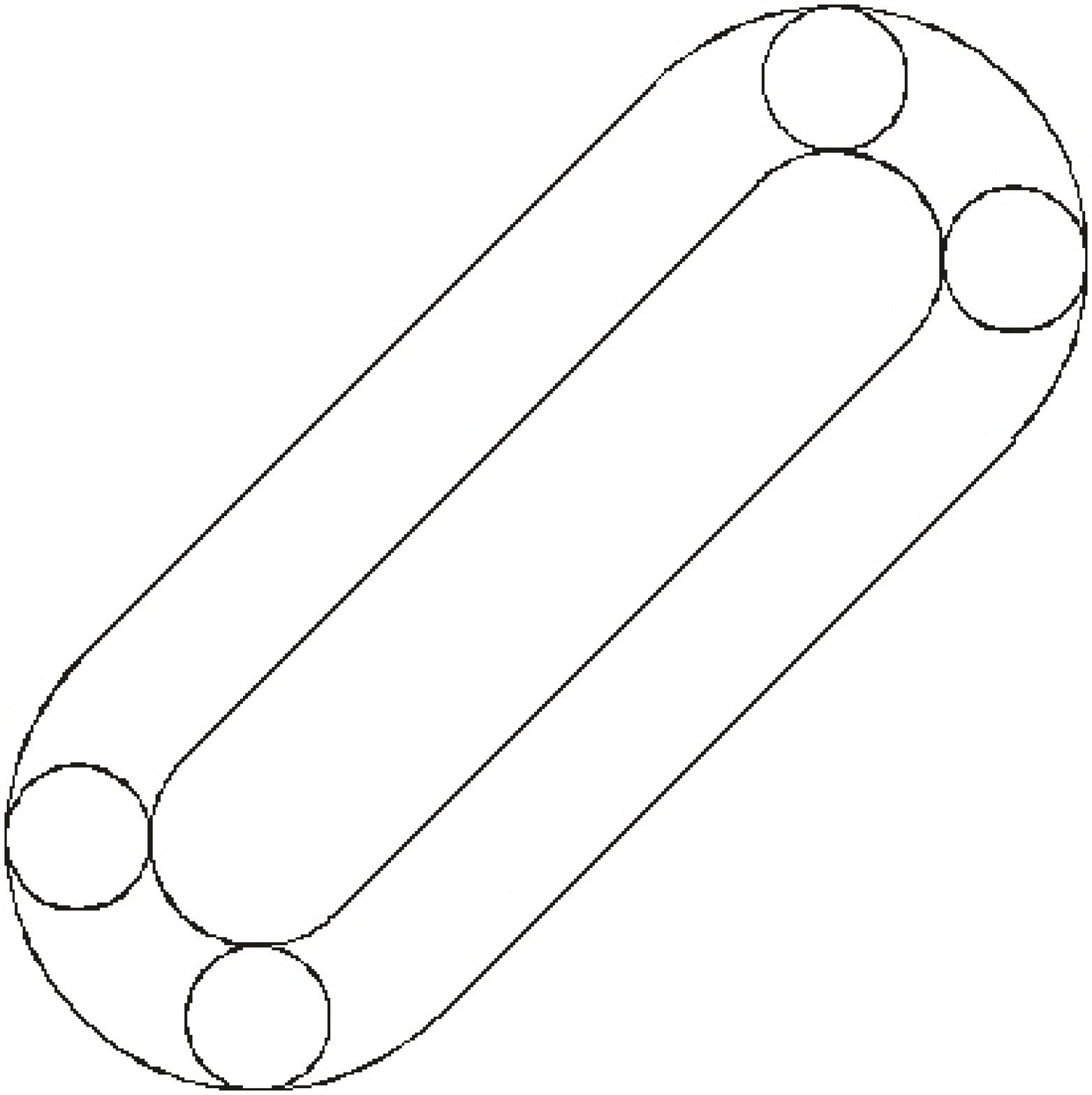}%
}%
{\includegraphics[
width=1.2503in
]%
{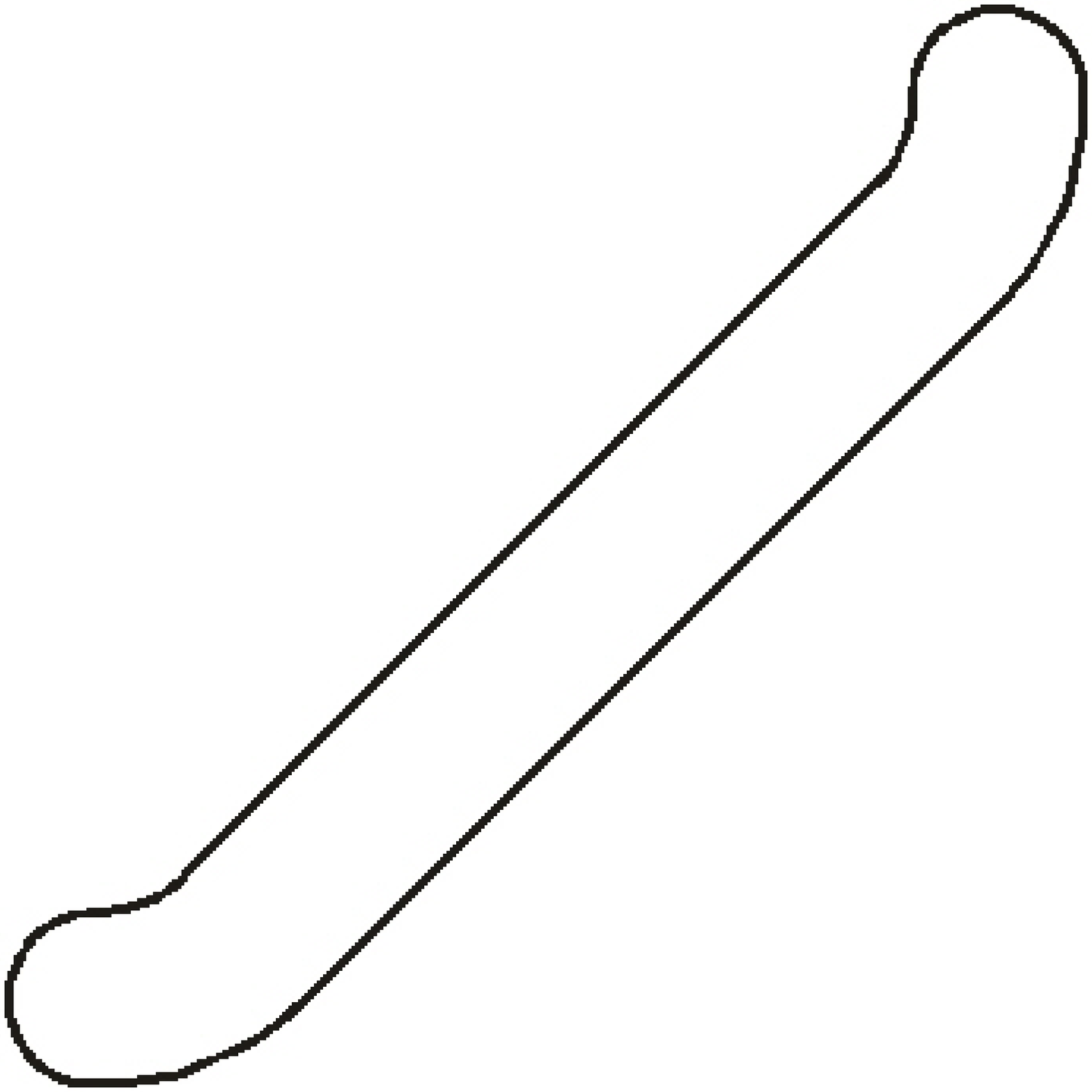}%
}%
%

\begin{center}
\includegraphics[
width=1.2926in
]%
{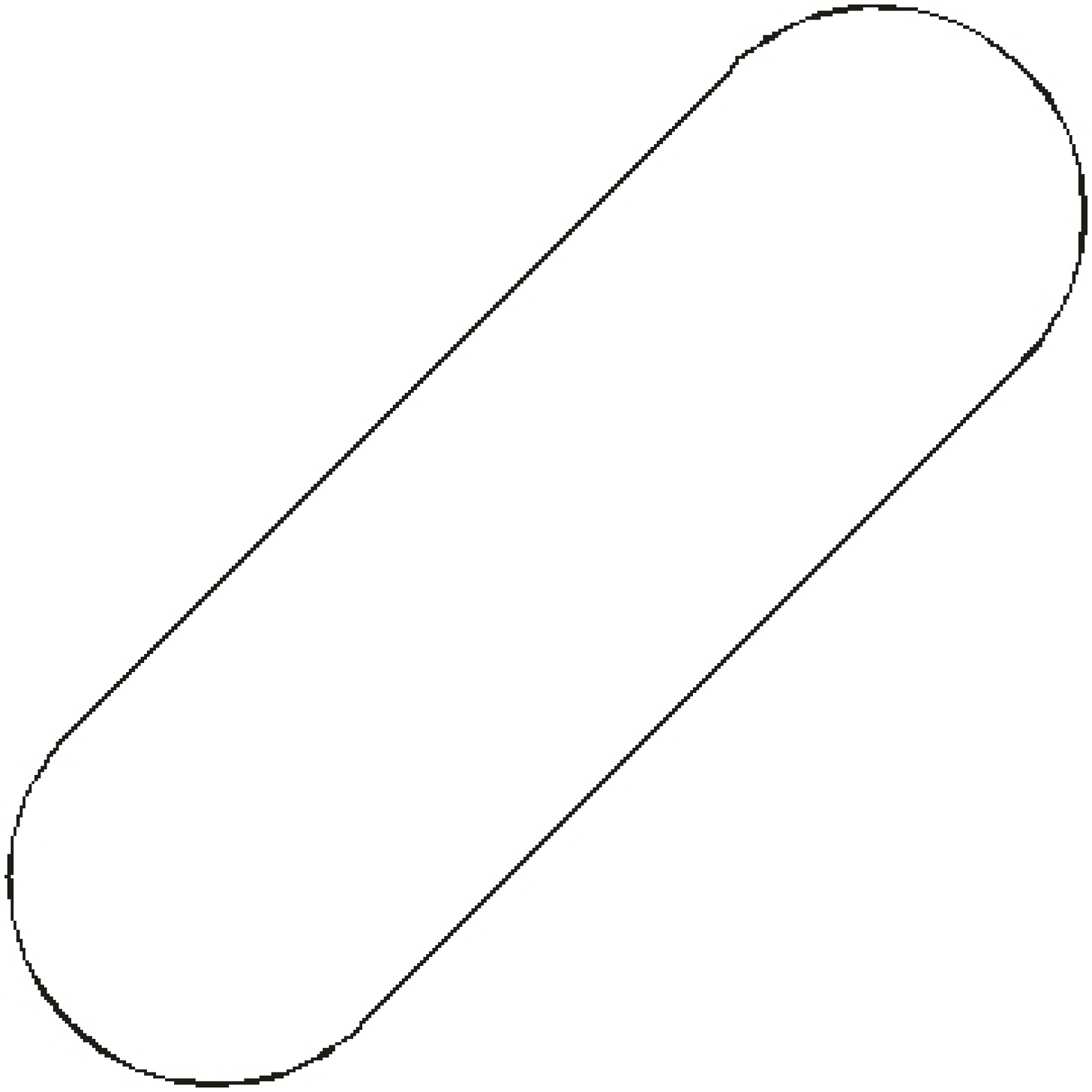}%
\end{center}

\textbf{Figure 8: }$\mathbf{M\in}\mathbb{R}^{2}$\textbf{, }$\widetilde
{\mathbf{M}_{\mathbf{F}}},$\textbf{\ }$\widetilde{\mathbf{M}}_{\mathbf{F}}%
\cup\left(  \underset{i=1}{\overset{n}{\cup}}\mathbb{S}_{i}^{1}\right)
=\underset{j=1}{\overset{m}{\cup}}\mathbb{S}_{j}^{1}$ and an $\mathbb{S}%
_{j}^{1}$.
\end{center}

The approach we take to solve both problems has two steps. The first step
eliminates cancellation due to two submanifolds in $\mathbb{R}^{n}$ with
opposite orientation. We endow the double cover with a canonical orientation.
So for example if a submanifold is a M\"{o}bius band, its double cover union
boundary lift will be a torus. The orientation on the lift is not inherited
from the underlying submanifold in $\mathbb{R}^{n}$. The second step is
eliminate cancellation due to degree zero Gauss maps under homothetic
contraction to a point of the submanifold in $\mathbb{R}^{n}$. This means we
must eliminate tori. We do this by adding topological cylinders so that the
result is a union of $\mathbb{S}^{n-1}$.

\begin{center}%
\[
\widetilde{\mathbf{M}}_{\mathbf{F}}\cup\left(  \underset{i=1}{\overset{n}%
{\cup}}\mathbb{S}_{i}^{n-1}\right)  =\underset{j=1}{\overset{m}{\cup}%
}\mathbb{S}_{j}^{n-1}%
\]

\end{center}

The induced orientation from step 1 applied to the extra parts of the lift
will have the orientation on each $\mathbb{S}^{n-1}\mathbb{\ }$that will
ensure a degree one Gauss map under homothetic contraction. Thus both mass
elimination problems are solved. See figure 8.

We need to check that the extra mass added with the union of $\mathbb{S}%
^{n-1}s$ is uniformly bounded for a subsequence of lifts.$\widetilde
{\mathbf{M}}_{\mathbf{Fi}}$. We need to start with a lemma:

\begin{lemma}
There exists a polyhedral n-2 skeleton S$_{k}$ of parts of great spheres on
$\mathbb{S}^{n-1}$ such that its pull-back to the lifts, $\widetilde{M}_{i}s$
(under the projection map g:$\mathbb{R}^{n}$\textsf{X}$\mathbb{S}%
^{n-1}\rightarrow\mathbb{S}^{n-1}),$ will be of uniformly bounded n-2
Hausdorff measure, and that the diameter of each n-1 cell in $\mathbb{S}%
^{n-1}$-S$_{k}$ is less than $1$.
\end{lemma}

\begin{proof}
This proof uses the coarea formula to express the ($n-1$) volume of an
$\widetilde{M}_{i}$ in terms of the ($n-2$) volume of co-dimension one sets on
$\widetilde{M}_{i}s.$ The codimension one sets form a disjoint 1-parameter
family that sweeps out to cover part of the volume of $\widetilde{M}_{i}$. In
this case the choice of 1-parameter family and the fact that $\widetilde
{M}_{i}$ is a lift ensures that the sweep rate is at least one. The volume of
$\widetilde{M}_{i}$ is finite. Therefore the $(n-2)$ volume of the codimension
one sets will be uniformly bounded away from some set of measure
$\varepsilon>0$ in the 1-parameter family. Otherwise the swept out volume of
$\widetilde{M}_{i}$ would be too great. This proves that there must be
candidates for the positions of the arcs of the skeleton $S_{k} $.

Take a family, parameterized by $t$ of great spheres passing through the point
$p$ on the $\mathbb{S}^{n-1}.$ Each \ great sphere corresponds to a fixed $t$,
and is parameterized by $\mathbf{x}.$ In the center of figure 7 parts of two
such great spheres are shown from $g(a)$ to $g(d)$ and $g(b)$ to $g(c)$. Now
remove two neighborhoods of $p$ and its antipode from the family of spheres.
The family of spheres now lies in the region shown by $g(a)$, $g(b)$, $g(c)$,
and $g(d)$.%

\begin{center}
\includegraphics[
width=3.7243in
]%
{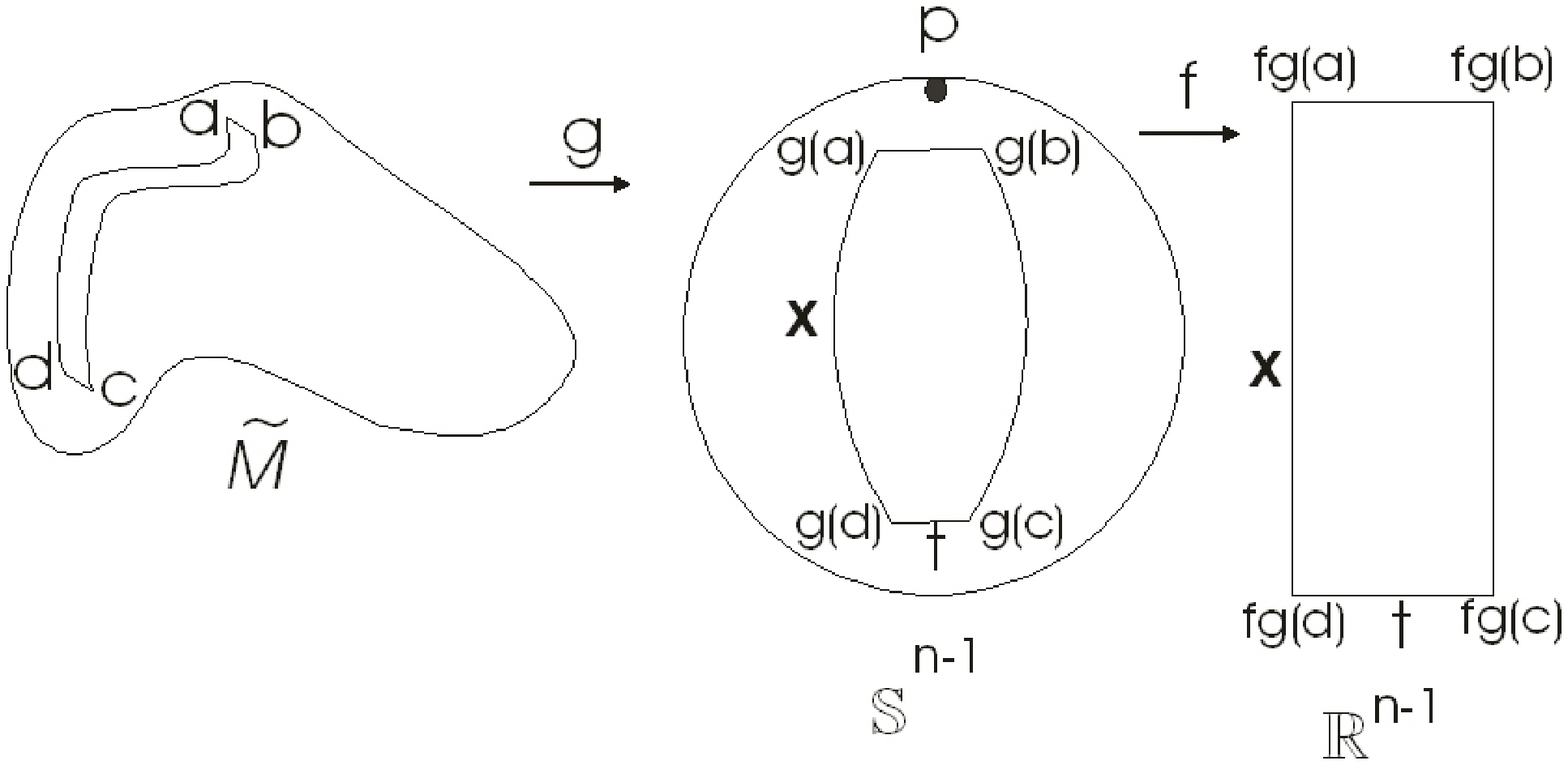}%
\\
Figure 9
\end{center}

Now apply the coarea formula on the projection map to the pre-image of the
great sphere sections parametrized by $t$ in $\widetilde{M}$. We need to map
from $\mathbb{R}^{n-1}$ to $\mathbb{R}$, this is projection $p$ into the $t$
coordinate. $(pfg)^{-1}(t)$ will be the pre-image in $\widetilde{M}$ of all
points with the given $t$ value over the values of $\mathbf{x}$. For example
the arcs $ad$ and $bc$ in $\widetilde{M}$ are each pre-images for the max and
min values of $t$. The coarea formula [MF: p30][S: p53][F1: 3.2] is:%

\begin{equation}
\underset{\widetilde{M}}{\int}J_{pfg}dH^{n-1}=\underset{t\in\lbrack0,1]}{\int
}H^{n-2}\left(  (pfg)^{-1}\left(  t\right)  \right)  dL^{1}\left(  t\right)
\end{equation}
\qquad

We now need to find J$_{g}$, J$_{f}$ and J$_{p}$ in the direction orthonormal
to the level sets of $t$. As g: $\widetilde{M}\subset\mathbb{R}^{n}$%
\textsf{X}$\mathbb{S}^{n-1}\rightarrow\mathbb{S}^{n-1}$, is a projection map,
$J_{g}$ $\leq1.$ We know by construction that $J_{f}\leq1,$ and $J_{p}=1$.

$J_{pfg}\leq1.$ So we can now write%

\begin{equation}
\left\vert \widetilde{M}\right\vert =\underset{\widetilde{M}}{\int}%
dH^{n-1}\geq\underset{t\in\lbrack0,1]}{\int}H^{n-2}\left(  (pfg)^{-1}\left(
t\right)  \right)  dL^{1}\left(  t\right)
\end{equation}
\qquad

We can conclude that $(pfg)^{-1}\left(  t\right)  $ is finite for almost all
$t $. Furthermore we know that for every bounded B there will be a set of
positive measure $e$ in [0,1] for values of $t$ where $H^{n-2}\left(
(pfg)^{-1}\left(  t\right)  \right)  <B.$ Let $H^{n-2}\left(  (pfg)^{-1}%
\left(  t\right)  \right)  \geq B$ on a set of measure s. Now we can write
$sB\leq\left\vert \widetilde{M}\right\vert ,$ and as $e=1-s,$ we can conclude
$e=1-\frac{\left\vert \widetilde{M}\right\vert }{B}.$

For any fixed $\widetilde{M}$ we can now construct the desired skeleton on
$\mathbb{S}^{n-1}$ by finding a union of arcs. Each arc is positioned
transversely according to a $t$ value for which it will have a pullback to
$\widetilde{M}$ of bounded volume.

We now need to show that we can derive a uniform bound for all $\widetilde
{M}_{i}$. For this we need to pass to a subsequence. We know that for each $i$
a portion of the interval [0,1] of finite measure, $H^{n-2}\left(
(pfg)^{-1}\left(  t\right)  \right)  <B$. Now we consider the sum of the
characteristic functions of these sets on [0,1]. The integral of this sum over
[0,1] must be infinite. We can argue by contradiction that if there is no good
direction then the sum will be finite over the whole of [0,1] which gives a
contradiction. Therefore there will be at least a set of positive measure on
[0,1] for values of $t$, where an infinite number of $i^{\text{'}} $s have
those $t$'s as good directions. We therefore can pass to this subsequence
$\widetilde{M}_{i_{j}}$ to gain the uniform bound as required.
\end{proof}

\begin{lemma}
The pull-back of the intersections of great spheres on $\mathbb{S}^{n-1}$ in
lemma 5.1 can also be of uniformly bounded mass for a subsequence of the
$\widetilde{M}_{i}$.

\begin{proof}
This follows by the same argument using the coarea formula and choice of
spheres. Although the intersections have lower dimension, there is a
correspondingly higher parameter family of them. So the dimensions work out in
the coarea formula. As there are a finite number of codimensions we only need
to take subsequences of the $\widetilde{M}_{i}$ a finite number of times. Thus
the process terminates to give a subsequence.
\end{proof}
\end{lemma}

\begin{lemma}
The extra mass of the spheres added is uniformly bounded by uniform linear
multiples of masses of the pull-backs in lemmas 5.1 and 5.2.
\end{lemma}

\begin{proof}
Consider a fiber in $\mathbb{R}^{n}$\textsf{X}$\mathbb{S}^{n-1}$ above a point
in $\mathbb{R}^{n}.$ Suppose it is on the pull-back of part of a great sphere
as in lemma 5.1. There are two antipodal arcs on that $\mathbb{S}^{n-1}$ which
will contain antipodal points in $\widetilde{M}_{F}.$ We can connect them in
$\mathbb{S}^{n-1}$ by the extending the arcs on which they lie until they met
in one direction. This will have a bounded arc length. As we move along the
pull-back of the S$_{k}$ from lemma 5.1, we always connect the same way. This
means we are connecting by adding part of a product structure, that is (the
pull-back)$\mathsf{X}$(an arc in S$^{n-1}$). This means that however the
pull-back may curve around, this curvature does not add to the mass of the
added S$^{n-1}$. Whatever the dimension of the fiber that is needed for
filling in, and whatever the dimension of the pull-back set of spheres or
their intersections we always use a product structure to ensure that the added
mass is fixed linear multiple of the mass of the pullbacks from lemmas 5.1 and 5.2.
\end{proof}

\textbf{We can now proceed to theorem 5.2.}

\begin{theorem}
When the subsets of $\mathbb{R}^{n}$ meet conditions 5.4.i, 5.4.ii and
5.4.iii, we can recover all collapsed dimensions in limit union of subsets of
$\mathbb{R}^{n}$
\end{theorem}

\begin{proof}
To ensure we do not lose any sets under dimension collapsing under homothety
we perform a surgery on the lift that makes it into a union of embedded
spheres as in figure 9. We need to check this can always be done without
adding more than a bounded amount of mass. Lemmas 5.1, 5.2 and 5.3 above do this.

These spheres then will have of degree of projection map of plus 1 or minus 1
depending on the orientation. We choose the orientation of the sphere so that
the degree is always plus 1. We prove this can be done by showing that the
projection map of an embedded sphere is isotopic to the identity. We can
continuously deform the sphere in $\mathbb{R}^{n}$\textsf{X}$\mathbb{S}^{n-1}
$ to make it convex, then contract down to a point.

Now we need to verify that this construction will give us a residual measure
for both full and partial dimensional collapsing. The full collapsing case is
given by the degree of projection map being 1. Now to prove the partial
dimension collapsing case, You get a $\mathbb{R}^{i}$\textsf{X}$\mathbb{S}%
^{j-i}$ as a local product structure in the limit. Now apply the same isotopy
argument in those dimensions which collapse under the product structure locally.

We finally need to check that this assignation of orientation does not induce
dimension collapsing with homothetic contraction. At each point in
$\mathbb{R}^{i}$\textsf{X}$\mathbb{S}^{n-1}$ the assigned orientation is
canonical for each tangent plane that is a lift or a fill in in the fibre of a
point. It is the canonical outward pointing vector to the $\mathbb{S}^{n-1} $ component.
\end{proof}

\section{Summary}

We have shown that under suitable uniform bounds (i.e.:

\textbf{Either}

(\textit{(3.3.i) Sets and their boundaries, and lower dimensional faces, are
}$C^{2}$\textit{\ rectifiable}.

\textit{(3.3.ii) Sets and their boundaries, and lower dimensional faces, have
uniformly bounded finite mass and finite integrals of principle curvatures
(including distributional)}

\textit{(3.3.iii) Sets have bounds on the integrals of the derivatives of
curvature.})

\textbf{or}

(\textit{(5.4.i) Sets and their boundaries, and lower dimensional faces, are
}$C^{2}$ \textit{rectifiable}.

\textit{(5.4.ii) Sets and their boundaries, and lower dimensional faces, have
uniformly bounded finite mass and finite integrals of principle
curvatures(including distributional).}

\textit{(5.4.iii) The topologically singular set, such as Y-singularities,
must have a uniformly finite mass lift.})),

a sequence of unions of C$^{2}$ j-rectifiable subsets of $\mathbb{R}^{n}$,
with dimensions $j$ ranging from 0 to n-1, will converge to a union of
rectifiable subsets of $\mathbb{R}^{n}$ with dimensions ranging from 0 to n-1,
under a topology based upon the outward pointing vectors. This topology will
also capture dimension collapsing.

Subsets whose lifts contain infinite $H^{n-2}$ measure stationary non-manifold
sets such as honeycombs in appendix may possibly converge only using varifold
convergence on $\mathbb{R}^{n}$\textsf{X}$\mathbb{S}^{n-1}$. Meanwhile others
subsets such as sequences of polyhedral approximations may have lifts that
possibly converge only with current compactness. When all subsets are
boundariless immersed manifolds, current compactness on the lifts is more
general, as it does not have hypotheses (3.3.iii) that place bounds on
derivatives of curvature.

\section{Appendix}

We give an example of a sequence of rectifiable sets with finite first
variation and finite mass, but infinite homological boundary.

We will describe first variation and homological boundary.

\subsection{First Variation}

The first variation of an m-submanifold $U$ in $\mathbb{R}^{n}$ associated
with a smooth compactly supported vector field $\mathbf{g}$ on $\mathbb{R}%
^{n}$ is:

$%
{\displaystyle\int\limits_{U}}
\mathbf{H}.\mathbf{g}dH^{m}+%
{\displaystyle\int\limits_{\partial U}}
-\mathbf{\nu}.\mathbf{g}dH^{m-1}$, is the total first variation where
$\mathbf{H}$ is the mean curvature vector on $U$, vector $\mathbf{\nu}$ is the
inward pointing normal vector on the boundary.

Notice that at a Y-singularity where three half planes come together the first
variation terms from each half plan can cancel.

\subsection{Homological Boundary}

This is the usual boundary in homology theory. It is the notion of boundary
suitable for currents. It coincides with set boundary as in the first
variation, but has extra terms due to orientation that may not match up. For
example a M\"{o}bius band can be represented as a current with boundary. There
will have to be an extra boundary component of multiplicity 2 due to the non
orientability of the M\"{o}bius band.

The Y-singularities will have one homological boundary component at least for
usual currents and homology. Taylor [T][MF: p105] uses currents mod n as a way
of getting around this problem.

\subsection{Example of varifold only convergence}

\begin{center}%
\begin{center}
\includegraphics[
width=1.2644in
]%
{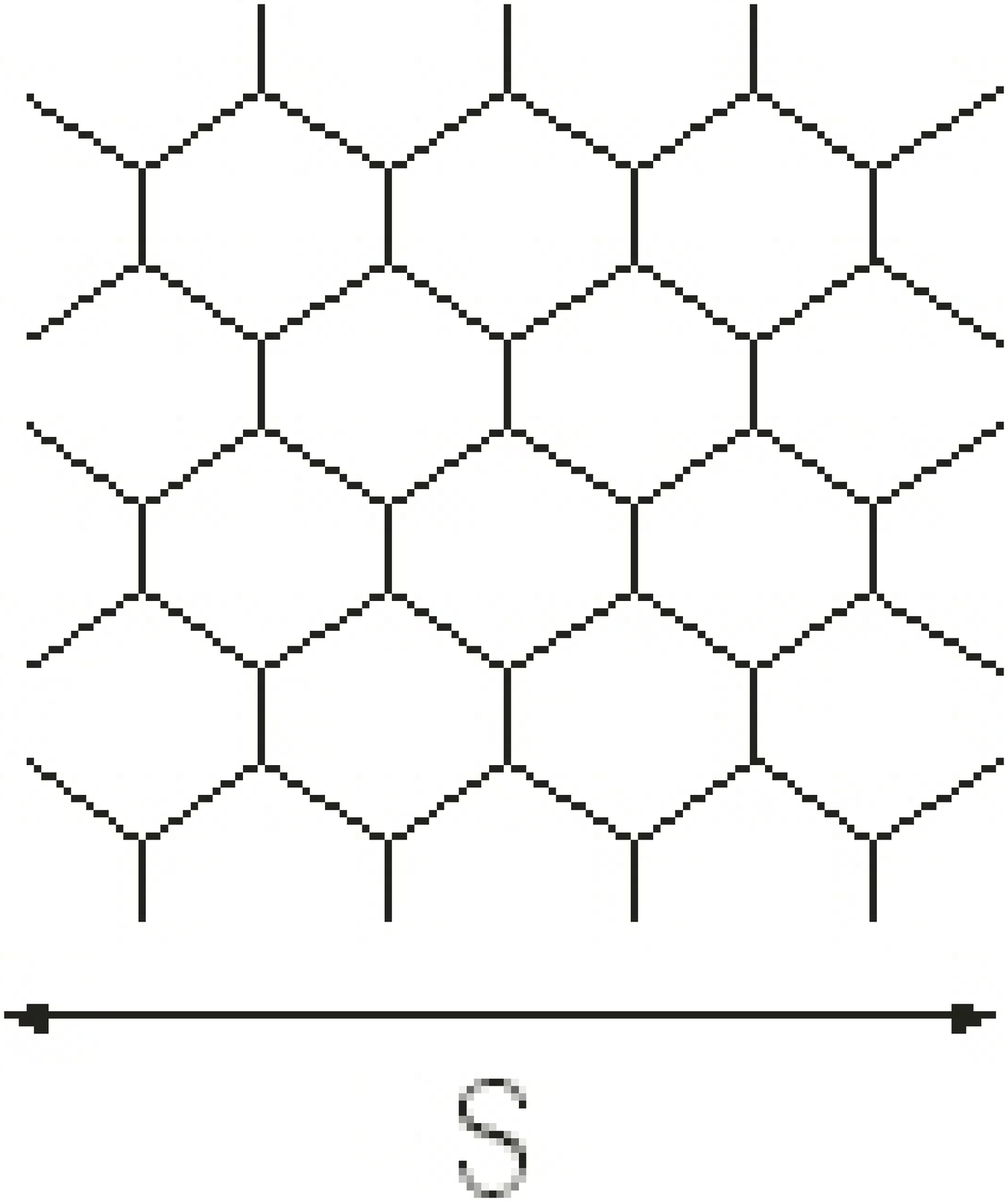}%
\end{center}

Figure 10
\end{center}

Consider a honeycomb, figure 10, in a cube of edge length $s$ in
$\mathbb{R}^{3}$. It projects down to a tessellation of hexagons in the $x-y$
plane and has height $h$ in the $z$ direction. Along each side of the
honeycomb there are $n$ vertical edges. This will correspond to the order of
$n^{2}$ vertical Y-singularities on the interior of the honeycomb, where three
faces come together.

The total mass of the honeycomb is of the order $2nsh$. As $s$ does not affect
any other quantity, we can now ignore mass.

The homological boundary mass, from the Y's is of the order of $hn^{2}$. The
first variation, from the outer ends of the honeycomb is of order $4nh$. We
can set up a sequence where $h=1/n$ of honeycombs. The first variation is
uniformly bounded and the homological boundary is not.

Notice that if $h=1/n^{3}$ we can take the unions of all the honeycombs and
obtain a rectifiable set which can be represented as an integer multiplicity
rectifiable varifold, but not as an integer multiplicity rectifiable current.

This particular example can be represented as a current mod 3, but one can
easily add extra faces to the honeycomb to create an infinite measure of
singularities with 3 and infinite measure with 5 faces coming together. Thus
in general currents mod n do not eliminate the need for varifold compactness.

\textbf{Acknowledgements}

I would like to thank Bob Hardt, my thesis advisor who suggested the approach
for this work, Thierry de Pauw, Robert Gulliver, Brian White, Penny Smith, and
David Johnson for their helpful discussions.

\end{document}